\documentclass[a4paper,12pt]{article}

\usepackage[margin=2.5cm]{geometry}
\usepackage{graphicx}
\usepackage{amssymb}

\usepackage{amsmath}
\usepackage{empheq}
\usepackage{slashed}
\usepackage{relsize}

\usepackage{float}
\usepackage{caption}
\usepackage{subcaption}
\usepackage{comment}
\usepackage{xcolor}

\usepackage[nosort]{cite}

\usepackage[colorlinks=true,
			bookmarks=true,
			citecolor=black,
			linkcolor=black,
			urlcolor=black,
			hypertexnames=false]{hyperref}

\newcommand{\dif}{\mathrm{d}}

\newcommand{\ksum}{%
	\mathop{\ooalign{$\sum$\cr\hidewidth$\mkern3mu \raisebox{.3ex}{\scriptsize $k$}\mkern3mu$\hidewidth}}%
}

\newcommand{\bigksum}{%
	\mathop{\ooalign{$\mathlarger{\sum}$\cr\hidewidth$\mkern3mu \raisebox{.3ex}{\scriptsize $k$}\mkern3mu$\hidewidth}}%
}

\newcommand{\kplusonesum}{%
	\mathop{\ooalign{$\mathlarger{\sum}$\cr\hidewidth$\mkern3mu \raisebox{.3ex}{\scriptsize $k+1$}\mkern3mu$\hidewidth}}%
}

\newcommand{\kminusonesum}{%
	\mathop{\ooalign{$\mathlarger{\sum}$\cr\hidewidth$\mkern3mu \raisebox{.3ex}{\scriptsize $k-1$}\mkern3mu$\hidewidth}}%
}

\newcommand{\fracb}[2]{\displaystyle{\frac{#1}{#2}}}

\colorlet{mycolor}{black}

\allowdisplaybreaks

\makeatletter
\renewcommand\section{\@startsection {section}{1}{\z@}%
	{-3.5ex \@plus -1ex \@minus -.2ex}%
	{2.3ex \@plus.2ex}%
	{\setcounter{equation}{0}\large\bfseries}}

\renewcommand\subsection{\@startsection{subsection}{2}{\z@}%
	{-3.25ex\@plus -1ex \@minus -.2ex}%
	{1.5ex \@plus .2ex}%
	{\normalfont\bfseries}}

\renewcommand\subsubsection{\@startsection{subsubsection}{3}{\z@}%
	{-3.25ex\@plus -1ex \@minus -.2ex}%
	{1.5ex \@plus .2ex}%
	{\normalfont\bfseries
}}
\makeatother

\makeatletter
\newcommand{\myheading}[1]{%
	\par\addvspace{2ex}                
	\noindent{\bfseries\normalsize #1}     
	\par\nobreak\smallskip            
	\@afterheading                    
}
\makeatother

\begin{document}

\thispagestyle{empty}
\begin{flushright}

\end{flushright}
\vbox{}

\vspace{2cm}

\begin{center}
{\LARGE{One-sided type-D Ricci-flat multi-centre metrics
}}\\[16mm]
{{Yu Chen}}
\\[6mm]
{\it School of Physical and Mathematical Sciences,
\\[1mm]
Nanyang Technological University, Singapore 637371}\\[15mm]

\end{center}
\vspace{1cm}

\centerline{\bf Abstract}
\bigskip
\noindent
{We employ a simplified form of Tod's ansatz---which is built on the LeBrun--Tod ansatz---to construct} one-sided type-D Ricci-flat multi-centre metrics. These metrics are Hermitian non-K{\"a}hler and conformally K{\"a}hler, and they are determined by a pair of generating potentials: an axisymmetric harmonic function {(different from the one originally proposed by Tod)} in an auxiliary 3D flat space and its ``harmonic conjugate''. We carry out a systematic study of these multi-centre metrics, including their rod structure, asymptotic structure and recovering from them various known closed-form examples. In particular, we show that, when fitted into the scheme of multi-soliton solutions on flat space constructed by the present author using the inverse-scattering method, these metrics, with number of centres $n\le 3$, have a free-soliton number $2$ or $3$ greater than their phantom-soliton number; we conjecture that this is true for general $n$.


\newpage

\tableofcontents

\section{Introduction}

Gravitational instantons are completely regular solutions to the vacuum Einstein field equations with Euclidean signature and no cosmological constant. They were originally studied in Euclidean quantum gravity and were proposed to make dominant contributions to the tunneling amplitudes between different vacua in the theory  \cite{Gibbons:1976ue,Hawking:1976jb}. They are also of interest in differential geometry as complete, non-compact Ricci-flat Riemannian four-manifolds, providing explicit geometrical models that illustrate deep structures like curvature decay and singularity resolution.

In a Riemannian manifold of dimension 4, the Weyl tensor splits into a self-dual part $W^+$ and an anti-self-dual part $W^-$ under the Hodge star operator $\star$, which is an involution. Each part/side of the Weyl tensor, $W^+$ and $W^-$, falls into one of the three Petrov types \cite{Przanowski:1983ybm,Derdzinski,Goldblatt:1994rx}: type I, type D and type O. When the (anti-)self-dual Weyl tensor, when viewed as a trace-free endomorphism on (anti-)self-dual two-forms, has three distinct eigenvalues, it is said to be of Petrov type I, i.e., the general type. When two among the three eigenvalues coincide, the (anti-)self-dual Weyl tensor is said to be of Petrov type D. For the remaining case with three identical and thus vanishing eigenvalues, the (anti-)self-dual Weyl tensor is said to be of Petrov type O. The Petrov types of the two parts/sides of the Weyl tensor are independent of each other. We can classify gravitational instantons by their Petrov type. A gravitational instanton is of Petrov type X$^+$Y$^-$, for some X, Y $\in\{\text{I, D, O}\}$, where X denotes the Petrov type of $W^+$ and Y denotes the Petrov type of $W^-$. With orientation reversed, the Petrov type X$^+$Y$^-$ of a gravitational instanton will change to Y$^+$X$^-$.

The class of hyper-K{\"a}hler gravitational instantons have been studied extensively over the last few decades and now a complete classification of it is available. With the assumption of a tri-holomorphic Killing vector, hyper-K{\"a}hler gravitational instantons can be explicitly constructed using the so-called Gibbons--Hawking ansatz \cite{Gibbons:1978tef}. The (multi-centre) Taub--NUT \cite{Newman:1963yy,Hawking:1976jb} and Gibbons--Hawking \cite{Gibbons:1978tef} instantons belong to this subclass, and they are asymptotically locally flat (ALF) and asymptotically locally Euclidean (ALE) respectively. With all centres collinearly aligned, which we shall assume from now on, the Taub--NUT and Gibbons--Hawking instantons have toric isometry. The Eguchi--Hanson instanton \cite{Eguchi:1978xp} is a special case of the Gibbons--Hawking instantons with two centres and equal charges. The hyper-K{\"a}hler condition implies Ricci-flatness and a vanishing self-dual Weyl tensor (with orientation appropriately chosen). So the Taub--NUT and Gibbons--Hawking instantons are anti-self-dual and of Petrov type O$^+$X$^-$, for some X which may depend on the number of centres and other parameters. They are one-sided type O.

There is another class of gravitational instantons whose metrics are known in closed form. The Schwarzschild\footnote{We shall work exclusively in Euclidean signature. We always mean the Euclidean section when we refer to a solution which has both Lorentzian and Euclidean sections.} instanton is an example in this class. Other examples in this class include the Kerr instanton \cite{Gibbons:1976ue}, the Taub-bolt instanton \cite{Page:1978hdy}, and the more recently discovered Chen--Teo instanton \cite{Chen:2011tc}. By the work of Aksteiner and Anderson \cite{Aksteiner:2021fae}, the Chen--Teo instanton is known to be one-sided type D. It is then clear that this class of gravitational instantons are all one-sided type D: they all have Petrov type D$^+$X$^-$ for some X (with orientation chosen appropriately). Biquard and Gauduchon \cite{Biquard:2021gwj} provided a complete classification of toric one-sided type-D (Hermitian non-K{\"a}hler) ALF gravitational instantons. With the inclusion of the single-centre Taub--NUT instanton with orientation opposite to the hyper-K{\"a}hler orientation, the above list is in fact complete.

From a constructive point of view, to find gravitational instantons, we often first find local metrics/solutions\footnote{The terms ``metrics'' and ``solutions'' (to Einstein's field equations) shall be used interchangeably in this paper.} that solve the vacuum Einstein field equations and then impose necessary regularity conditions, such as the absence of curvature and conical singularities. For example, in the hyper-K{\"a}hler (and thus one-sided type-O) case, the Taub--NUT and Gibbons--Hawking instantons are obtained from the multi-centre Gibbons--Hawking metrics \cite{Gibbons:1978tef}. When all the centres align on an axis with charges not all equal, the Gibbons--Hawking metrics contain mild singularities: there are conical singularities along the axis and also potential orbifold singularities at where adjacent segments of the axis meet. When all the charges are set equal, these mild singularities on the axis can be removed and we obtain gravitational instantons. In the one-sided type-D case, the Kerr and Taub-bolt instantons can be extracted from the Kerr--NUT metric \cite{Carter:1968ks,Gibbons:1979nf} by imposing absence of conical and orbifold singularities. The Chen--Teo instanton can be extracted from the Chen--Teo metric \cite{Chen:2015vva}.

The Taub--NUT and Gibbons--Hawking instantons are all special cases of the multi-centre Gibbons--Hawking metrics, while one-sided type-D gravitational instantons were obtained from different families of metrics. This motivates us to ask the following question: What is the one-sided type-D counterpart of the multi-centre Gibbons--Hawking metrics? The field equations of one-sided type-D Ricci-flat metrics were reduced to a single second-order, nonlinear differential equation of 
one real function by Przanowski \cite{Przanowski:1984} and more recently by Tod \cite{Tod:2020ual}, and a symmetry was shown to exist. The form of metric that Tod obtained is known in the literature as LeBrun--Tod ansatz, which can be adapted from LeBrun's ansatz \cite{Lebrun:1991} for scalar-flat K{\"a}hler surfaces. With the assumption of an additional symmetry, {Tod \cite{Tod:2020ual} linearised the field equations by using a Ward transform and further reduced them to a single axisymmetric harmonic potential in an auxiliary 3D flat space. His final form of metric is known as Tod's ansatz.} In this paper, we {simplify the form of Tod's ansatz and} reduce the field equations to a pair of axisymmetric potentials, with one harmonic on the auxiliary 3D flat space and the other being the ``harmonic conjugate''.  {Our harmonic potentials are different from Tod's one and they are derivatives of the latter.} By choosing a simple form of these generating harmonic potentials, we obtain one-sided type-D Ricci-flat multi-centre metrics, which can be regarded as the counterpart of the multi-centre Gibbons--Hawking metrics {with all centres collinearly aligned.} 

We remark that one-sided type-D Ricci-flat multi-centre metrics are Hermitian and non-K{\"a}hler, and they appeared in equivalent forms and were approached from a different angle in earlier works \cite{Biquard:2021gwj,Biquard:2023rbr} of Biquard and Gauduchon with a primary goal to classify Hermitian non-K{\"a}hler ALF gravitational instantons. In this paper, our focus is on one-sided type-D Ricci-flat multi-centre metrics themselves. We perform a systematic study of these metrics, including, in particular, their construction, local geometric properties, and classification.

The paper is organised as follows. In Sec.~\ref{sec_type_D_metrics}, we first show that one-sided type-D Ricci-flat metrics admit a Killing vector and derive the LeBrun--Tod ansatz. Then by assuming an additional symmetry, {we derive Tod's ansatz in a simpler form in terms of a pair of generating potentials.} In Sec.~\ref{sec_multi_soliton_solutions}, we present one-sided type-D Ricci-flat multi-centre solutions. Their rod structures are studied, and it is then shown that they are ALF and admit two ALE limits. In the following section, we study some examples in detail and show that one-sided type-D Ricci-flat $n$-centre solutions for $n=1$, 2 and 3 are respectively the single-centre Taub--NUT solution, the Kerr--NUT solution and the Chen--Teo solution. In Sec.~\ref{sec_ISM_classification}, we discuss how these multi-centre solutions fit into the scheme of multi-soliton solutions on flat space constructed by the present author using the inverse-scattering method. In the last section, we discuss some potential extensions of the present work. The paper ends with an appendix, in which we give a simplified proof of a key result in Biquard--Gauduchon classification of Hermitian non-K{\"a}hler ALF gravitational instantons.

\section{{Toric} one-sided type-D Ricci-flat metrics}
\label{sec_type_D_metrics}

{In this section, we give a self-contained derivation of Tod’s ansatz of toric one-sided type-D Ricci-flat metric. We present the metric in a form slightly simpler than Tod's and with an alternative generating potential.} 

\subsection{The LeBrun--Tod ansatz}

{In \cite{Tod:2020ual}, Tod first derived an ansatz, the Lebrun--Tod ansatz, for one-sided type-D Ricci-flat metrics (not necessarily toric).} In this subsection, we follow Tod \cite{Tod:2020ual,Tod:2024msz} to rederive the LeBrun--Tod ansatz using the  two-component spinor formalism. We follow Penrose and Rindler \cite{Penrose:1985bww,Penrose:1986ca} on basic notations of two-component spinors. For the adaption to Euclidean signature, we have referred to \cite{Dunajski:2010zz,Goldblatt:1994rx}.

\myheading{One-sided type-D condition}
Consider a Riemannian manifold in dimension four. We assume that the self-dual Weyl tensor $W^+$ of the metric has Petrov type D, so its spinorial image $\Psi_{ABCD}$ can be factorised as
\begin{equation}
	\label{self-dual-Weyl-tensor}
	\Psi_{ABCD}=6\Psi_2 o_{(A}o_B\iota_C\iota_{D)},
\end{equation}
where indices in round brackets are symmetrised. $\Psi_2$ is the only non-vanishing unprimed Weyl scalar. $(o^A,\iota^A=o^\dagger{}^A)$ form a spin frame with normalisation $o_A\iota^A=1$. $o^\dagger{}^A$ denotes the conjugate spinor of $o^A$ and $\iota^\dagger{}^A=-o^A$. Denote the primed spin frame by $(o^{A'}, \iota^{A'})$ with normalisation $o_{A'}\iota^{A'}=1$. A null tetrad $(l^a,n^a,m^a,p^a)$ is defined by the following correspondence
\begin{equation}
	l^a=o^Ao^{A'},\quad n^a=\iota^A\iota^{A'},\quad m^{a}=o^A\iota^{A'},\quad p^a=\iota^Ao^{A'},
\end{equation}
with normalisation
\begin{equation}
 l_an^a=1,\quad m_ap^a=-1,
\end{equation}
and complex conjugation
\begin{equation}
	\bar{l}^a=n^a,\quad \bar{m}^a=-p^a.
\end{equation}
The metric tensor $g_{ab}$ in terms of the null tetrad is given by
\begin{equation}
	g_{ab}=2l_{(a}n_{b)}-2m_{(a}p_{b)}.
\end{equation}

With the assumption of Ricci-flatness, by the Goldberg--Sachs theorem adapted to Euclidean signature (see, e.g., \cite{Goldblatt:1994rx}), the principal spinor $o_A$ of the Weyl spinor $\Psi_{ABCD}$ is geodesic and shear-free:
\begin{equation}
	\label{gsf_condition}
	o^Ao^B\nabla_{AA'}o_B=0\quad\Leftrightarrow\quad \kappa=\sigma=0.
\end{equation}
Here, $\kappa$ and $\sigma$ are two spin coefficients. Together with another two spin coefficients $\rho$ and $\tau$ that we will use later, they are defined as
\begin{equation}
	\begin{aligned}
		\kappa
		&=o^Ao^{A'}o^B\nabla_{AA'}o_B,\quad \sigma=o^A\iota^{A'}o^B\nabla_{AA'}o_B,\\
		\rho&=\iota^Ao^{A'}o^B\nabla_{AA'}o_B,\quad 	\tau=\iota^A\iota^{A'}o^B\nabla_{AA'}o_B.
	\end{aligned}
\end{equation}
The conjugate spinor $\iota^A$ is also geodesic and shear-free. The Bianchi identity in this case is $\nabla^A_{B'}\Psi_{ABCD}=0$, which gives a differential equation for $\Psi_2$:
\begin{equation}
	\label{Psi2_equation}
	\dif \ln \Psi_2=3(\bar{\rho}l+\rho n+\bar{\tau}m-\tau p),
\end{equation}
where a bar denotes complex conjugation as before.

We define an almost complex structure $J$ by its action on vectors $v^a$ as $(Jv)^b=J_{a}{}^bv^a$ with $J_a{}^b$ given by
\begin{equation}
	J_a{}^b=i\left(o_A\iota^B+\iota_Ao^B\right)\epsilon_{A'}{}^{B'},
\end{equation}
where $\epsilon_{A'B'}$ is the ``metric tensor'' in the space of primed spinors. Clearly, $J$ is real and satisfies $J^2=-\text{id}$. It can be checked that $J$ preserves the metric tensor $g_{ab}$. The Nijenhuis tensor $N_{bc}{}^a$ corresponding to $J$ can be calculated directly and can be written in terms of two-forms as
\begin{equation}
	N^a=4(\kappa n^a-\sigma p^a)p\wedge n+4(\bar{\kappa}l^a+\bar{\sigma}m^a)l\wedge m,
\end{equation}
which vanishes in view of (\ref{gsf_condition}). The almost complex structure $J$ is thus integrable. One-sided type-D Ricci-flat metrics are thus Hermitian. The fundamental two-form $\omega_{ab}$ defined by $\omega_{ab}=J_a{}^cg_{cb}$ is
\begin{equation}
	\label{fundamental_two_form}
	\omega_{ab}=2io_{(A}\iota_{B)}\epsilon_{A'B'}=[i(l\wedge n-m\wedge p)]_{ab}.
\end{equation}
It is not closed. The Lee form $\theta$ is defined to measure the non-closure of the fundamental two-form $\omega$:
\begin{equation}
	\label{closedness_Kahler_form2}
	\dif \omega =\theta \wedge\omega, 
\end{equation}
and it can be calculated as
\begin{equation}
\theta=-2(\bar{\rho}l+\rho n+\bar{\tau}m-\tau p)=-\frac{2}{3}\dif \ln \Psi_2.
\end{equation}
It follows that 
\begin{equation}
	\label{closedness_Kahler_form}
	\dif (\Psi_2^{\frac{2}{3}}\omega)=0.
\end{equation}
A one-sided type-D Ricci-flat metric $g_{ab}$ is thus Hermitian and non-K{\"a}hler, but conformally K{\"a}hler. It is conformal to the K{\"a}hler metric $\hat{g}_{ab}=\Psi_2^{\frac{2}{3}}g_{ab}$ with K{\"a}hler form $\hat{\omega}_{ab}=\Psi_2^{\frac{2}{3}}\omega$.

We have assumed that the metric is non-trivially one-sided type D, so $\Psi_2$ is not the zero function. A valence-2 symmetric Killing spinor $K_{AB}$ can be constructed (up to a multiplicative constant)
\begin{equation}
	K_{AB}\propto 2i\Psi_2^{-\frac{1}{3}}o_{(A}\iota_{B)},
\end{equation}
and it satisfies the equation
\begin{equation}
	\nabla_{C'(A}K_{BC)}=0.
\end{equation}
The Killing spinor $K_{AB}$ gives rise to a Killing vector $k^a$ defined by
\begin{equation}
	k^{AA'}=\frac{1}{3}\nabla_{B}^{A'}K^{BA}.
\end{equation}
The Killing one-form $k_a$ can be calculated as
\begin{equation}
	\label{Killing_vector_expression}
	k_a\propto i\Psi_2^{-\frac{1}{3}}(-\bar{\rho}l_a+\rho n_a-\bar{\tau}m_a-\tau p_a).
\end{equation}
It then follows that
\begin{equation}
	\label{J_on_k}
	k\lrcorner \omega\propto d\left(\Psi_2^{-\frac{1}{3}}\right),
\end{equation}
where $\omega$ is the fundamental two-form (\ref{fundamental_two_form}). The Killing vector $k^a$ is thus Hamiltonian.

\myheading{The LeBrun--Tod ansatz}

With the above results, we can set up coordinates, adapted to one-sided type-D Ricci-flat metrics. Define a function $W$ by $W^{-1}=k_ak^a$ and two normalised orthogonal one-forms $\theta^{(0)}$ and $\theta^{(1)}$ as
\begin{equation}
\label{theta_0and1}
\theta^{(0)}=W^{\frac{1}{2}}k,\quad	\theta^{(1)}=J\theta^{(0)}.
\end{equation}
In view of (\ref{J_on_k}), we can introduce a coordinate $w$ by
\begin{equation}
	\label{theta_1}
	\theta^{(1)}=W^{\frac{1}{2}}\dif w, \quad w\propto \Psi_2^{-\frac{1}{3}}.
\end{equation}
We now use a Lorentz transformation consisting of an $SU(2)$ transformation that mixes only primed spinors $(o^{A'},\iota^{A'})$ to set the spin coefficient $\tau=0$.  In view of (\ref{theta_0and1}), (\ref{Killing_vector_expression}), (\ref{theta_1}) and (\ref{Psi2_equation}), the complex span of $(\theta^{(0)},\theta^{(1)})$ is the same as that of $(l,n)$, and the complementary is spanned by $(m,p)$. The null-tetrad connection one-form
\begin{equation}
	-m_b\nabla_a l^b=-\tau l_a-\kappa n_a+\rho m_a+\sigma p_a=\rho m_a,
\end{equation}
is surface-forming, thanks to a result of Kerr (cf. Lemma (7.3.44) of \cite{Penrose:1986ca}) adapted to Euclidean signature. It follows that $m\propto P\dif \zeta$ for some complex functions $P$ and $\zeta$. By applying a Lorentz transformation consisting of an $SU(2)$ that rescales primed spinors by pure phases, we can make $P$ real. Therefore, we can write
\begin{equation}
	m=W^{\frac{1}{2}}e^{u/2}\dif \zeta,
\end{equation}
for some real function $u$.  We introduce two real coordinates $x$ and $y$ by $\zeta=x+iy$, and choose another two normalised  orthogonal one-forms $\theta^{(2)}$ and $\theta^{(3)}$, complementary to $\theta^{(0)}$ and $\theta^{(1)}$, as
\begin{equation}
	\theta^{(2)}=W^{\frac{1}{2}}e^{u/2}\dif x,\quad \theta^{(3)}=W^{\frac{1}{2}}e^{u/2}\dif y.
\end{equation}
Together, the one-forms $(\theta^{(0)},\theta^{(1)},\theta^{(2)},\theta^{(3)})$ form a real orthonormal tetrad. Lastly, we introduce a coordinate $\tau$ as the affine parameter along the flow generated by the Killing vector $k^a$, and thus
\begin{equation}
	\label{coordinate_tau}
	k^a\partial_a=\partial_\tau.
\end{equation}
From now on, the symbol $\tau$ shall be used as the coordinate defined in (\ref{coordinate_tau}), which bears no relation to the spin coefficient $\tau$. In the coordinate system ($\tau,w,x,y$) that we have introduced, the frame one-form $\theta^{(0)}$ can be written as
\begin{equation}
	\theta^{(0)}=W^{-\frac{1}{2}}(\dif \tau+A),
\end{equation}
for some one-form $A$ that lives in the quotient space of the four-manifold under study modulo the flow generated by the Killing vector $\partial_\tau$. The metric of the Riemannian four-manifold is then obtained as
\begin{equation}
	\label{metric_Lebrun_Tod}
    \dif s^2=\sum\limits_{s=0}^{3}\theta^{(s)}\theta^{(s)}=W^{-1}(\dif \tau+A)^2+W[\dif w^2+e^u(\dif x^2+\dif y^2)].
\end{equation}

With the above metric, now we impose the full one-sided type-D condition and Ricci-flat field equations. From the above construction of the orthonormal tetrad, we can see that the fundamental two-form (\ref{fundamental_two_form}) becomes
\begin{equation}
	\label{Kahler_form1}
	\omega=\theta^{(0)}\wedge\theta^{(1)}+\theta^{(2)}\wedge\theta^{(3)}=(\dif \tau+A)\wedge\dif w+e^u W\dif x\wedge\dif y.
\end{equation}
The distribution spanned by the (1,0)-forms
\begin{equation}
	\begin{cases}
		\theta^{(0)}+i\theta^{(1)}&=W^{-\frac{1}{2}}(\dif \tau+A)+iW^{\frac{1}{2}}\dif w\sim \dif \tau+A+iW\dif w,\\
		\theta^{(2)}+i\theta^{(3)}&=W^{\frac{1}{2}}e^{\frac{u}{2}}(\dif x+i\dif y)\sim \dif x+i\dif y,
	\end{cases}
\end{equation}
is integrable since $J$ is integrable. Expressions on two sides of $\sim$ differ by a multiplicative function, which is insignificant to generate the distribution. The integrability of the distribution leads to the following differential equation for $A$
\begin{equation}
	\label{Equation_for_A1}
	\dif A=-W_x\dif y\wedge\dif w-W_y\dif w\wedge\dif x-\eta \dif x\wedge \dif y,
\end{equation}
for some function $\eta$.
Eq.~(\ref{closedness_Kahler_form}) implies that the two-form $w^{-2}\omega$ is closed:
\begin{equation}
	\dif (w^{-2}\omega)=0.
\end{equation}
By substituting (\ref{Kahler_form1}) and (\ref{Equation_for_A1}) into this equation, we obtain
the complete equation for the one-form $A$ as
\begin{equation}
	\label{Equation_for_A2}
	\dif A=-W_x \dif y\wedge \dif w-W_y\dif w\wedge\dif x-w^2(w^{-2}e^uW)_w\dif x\wedge \dif y.
\end{equation}
The integrability of this equation is known as the monopole equation
\begin{equation}
	\label{monopolo_equation}
	W_{xx}+W_{yy}+[w^2(w^{-2}e^uW)_w]_w=0,
\end{equation}
a second-order differential equation for $W$. Thus, we see that for one-sided type-D Ricci-flat metrics, the one-form $A$ is completely determined by the function $W$  (up to an insignificant exact one-form).

In four dimensions, the space of two-forms splits into self-dual and anti-self-dual parts, each of which is three-dimensional. We can compute connection and curvature forms for these two parts separately. We choose orthonormal bases  $Z^{i}$ for the self-dual two-forms and $\tilde{Z}^{i'}$ for the anti-self-dual two-forms as
\begin{align}
	Z^{i}&=\frac{1}{\sqrt{2}}\left(\theta^{(0)}\wedge\theta^{(i)}+\frac{1}{2}\epsilon^i{}_{jk}\theta^{(j)}\wedge\theta^{(k)}\right),\label{self_dual_twoform_basis}\\
	\tilde{Z}^{i'}&=\frac{1}{\sqrt{2}}\left(\theta^{(0)}\wedge\theta^{(i)}-\frac{1}{2}\epsilon^i{}_{jk}\theta^{(j)}\wedge\theta^{(k)}\right).\label{antiself_dual_twoform_basis}
\end{align}
The symbol $\epsilon_{ijk}$ is completely anti-symmetric with respect to its indices and $\epsilon_{123}=1$ by definition. Since orthonormal bases of self-dual and anti-self-dual two-forms have been used, indices with respect to these bases are raised and lowered using Kronecker delta. The inner product of two two-forms $F$ and $G$ is defined as $\langle F,G\rangle\equiv \frac{1}{2}F_{ab}G^{ab}$, which has been used to normalise the basis two-forms $Z^{i}$ and $\tilde{Z}^{i'}$. 

We shall focus on the self-dual side only. The connection one-forms $\omega^{ij}$, defined by Cartan's first structure equation
\begin{equation}
	\dif Z^{i}+\omega^{i}{}_{j}{\wedge} Z^{j}=0,
\end{equation}
can be computed using the formula
\begin{equation}
	\omega^{ij}=\left(\star\dif Z^{i}\right)^\sharp \lrcorner Z^{j}-\left(\star\dif Z^{j}\right)^\sharp \lrcorner Z^{i}+\frac{1}{\sqrt{2}}\star\dif Z^{k}\epsilon_{k}{}^{{ij}}.
	\end{equation}
 Direct calculation using the above formula yields
\begin{align}
	\omega^{12}&=-w^{-1}W^{-\frac{1}{2}}\theta^{(2)},\\
	\omega^{23}&=W^{-\frac{1}{2}}\left[-w^{-1}\left(1-\frac{1}{2}wu_w\right)\theta^{(0)}+\frac{1}{2}u_ye^{-\frac{u}{2}}\theta^{(2)}-\frac{1}{2}u_xe^{-\frac{u}{2}}\theta^{(3)}\right],\\
	\omega^{31}&=w^{-1}W^{-\frac{1}{2}}\theta^{(3)}.
\end{align}
The curvature two-forms $F^{ij}$, defined by Cartan's second structure 
equation
\begin{equation}
	F^{ij}=\dif\omega^{ij}+\omega^{i}{}_{k}\wedge\omega^{kj},
\end{equation}
are found as
\begin{align}
F^{12}=&\sqrt{2}w^{-2}W^{-1}\left(1-\frac{1}{2}wu_w\right)Z^{3},\\
F^{23}=&-\frac{W^{-1}}{2\sqrt{2}}\left[e^{-u}\left((e^u)_{ww}+u_{xx}+u_{yy}\right)+8w^{-2}\left(1-\frac{1}{2}wu_w\right)\right]Z^{1}\nonumber\\
&-\frac{W^{-1}}{2\sqrt{2}}\left[e^{-u}(e^uu_{ww}-u_{xx}-u_{yy})-(2w^{-1}-u_w)(2w^{-1}-u_w-2W^{-1}W_w)\right]\tilde{Z}^{1'}\nonumber\\
&+\sqrt{2}e^{-\frac{u}{2}}\left[w^{-1}W^{-1}\left(1-\frac{1}{2}wu_w\right)\right]_x\tilde{Z}^{2'}\nonumber\\
&+\sqrt{2}e^{-\frac{u}{2}}\left[w^{-1}W^{-1}\left(1-\frac{1}{2}wu_w\right)\right]_y\tilde{Z}^{3'},\\
F^{31}=&\sqrt{2}w^{-2}W^{-1}\left(1-\frac{1}{2}wu_w\right)Z^{2}.
\end{align}

From the curvature two-forms $F^{ij}$, we can extract the self-dual Weyl tensor, the Ricci scalar and the trace-free Ricci tensor. The Ricci scalar is calculated as
\begin{equation}
	R=-e^{-u}W^{-1}\left[(e^u)_{ww}+u_{xx}+u_{yy}\right],
\end{equation}
the vanishing of which leads to the $SU(\infty)$-Toda equation for $u$:
\begin{equation}
	\label{Toda_3D}
  (e^u)_{ww}+  u_{xx}+u_{yy}=0.
\end{equation}
With this equation imposed, the self-dual Weyl curvature tensor, regarded as an endomorphism $W^+: \phi_{ab}\rightarrow \frac{1}{2}W^{+}_{ab}{}^{cd}\phi_{cd}$ on self-dual two-forms, takes the following form  in the basis (\ref{self_dual_twoform_basis}):
\begin{equation}
	\label{Weyl_tensor_self_dual}
	W^{+ij}=\begin{pmatrix}
		-2\Psi_2&0&0\\
		0&\Psi_2&0\\
		0&0&\Psi_2
	\end{pmatrix} \,\,\,\text{ where  }\,\,\,\Psi_2=w^{-2}W^{-1}\left(1-\frac{1}{2}wu_w\right),
\end{equation}
which is seen to be of Petrov type D.

With the function $u$ satisfying the $SU(\infty)$-Toda equation (\ref{Toda_3D}), the function $W$ satisfying the monopole equation (\ref{monopolo_equation}) and the one-form $A$ given by (\ref{Equation_for_A2}), the metric (\ref{metric_Lebrun_Tod}) with the fundamental two-form (\ref{Kahler_form1}) is scalar-flat and conformally K{\"a}hler. It can be adapted from LeBrun's ansatz \cite{Lebrun:1991} for scalar-flat K{\"a}hler surfaces.

With the scalar-flat condition (\ref{Toda_3D}) imposed, the trace-free Ricci curvature operator, from which the Ricci tensor can be reconstructed, takes the following form in the bases (\ref{self_dual_twoform_basis}) and (\ref{antiself_dual_twoform_basis}):
\begin{equation}
	\label{trace_free_Ricci}
	\Phi^{i'j}=\begin{pmatrix}
		w^{-2}\left[wW^{-1}\left(1-\frac{1}{2}wu_w\right)\right]_w&0&0\\
		e^{-\frac{u}{2}}w^{-2}\left[wW^{-1}\left(1-\frac{1}{2}wu_w\right)\right]_x&0&0\\
		e^{-\frac{u}{2}}w^{-2}\left[wW^{-1}\left(1-\frac{1}{2}wu_w\right)\right]_y&0&0
	\end{pmatrix}.
\end{equation}
The Ricci-flat condition requires $\Phi^{i'j}$ to vanish, giving rise to
\begin{equation}
	\dif \left[wW^{-1}\left(1-\frac{1}{2}wu_w\right)\right]=0,
\end{equation}
which can be integrated directly to give
\begin{equation}
	\label{W_equation1}
	W=\frac{1}{k_1}\left(w-\frac{1}{2}w^2u_w\right),
\end{equation}
where $k_1$ is an integration constant.\footnote{{We shall use positive $k_1$ in Eq.~(\ref{W_equation1}) without loss of generality, and our $w$ will be negative (see Eq.~(\ref{w_of_multisoliton})) to ensure $W>0$. Our convention is different from that used in Biquard and Gauduchon \cite{Biquard:2021gwj}. In their Eq.~(28), if the constant $k$ is positive then $\xi$ is also positive in order to ensure their $V$ (the equivalent of our $W$) to be positive. So our $w$ is the equivalent of their $-\xi$ (apart from a factor of 2).}} It can be checked that the function $W$ satisfies the monopole equation (\ref{monopolo_equation}). The Weyl scalar $\Psi_2$ (see (\ref{Weyl_tensor_self_dual})), with $W$ given by (\ref{W_equation1}), becomes
\begin{equation}
	\Psi_2=k_1w^{-3},
\end{equation}
as required from our construction (see Eq.~(\ref{theta_1})).

One-sided type-D Ricci-flat field equations have now been reduced to a single function $u$ that satisfies the $SU(\infty)$-Toda equation. The LeBrun--Tod ansatz is given by the metric (\ref{metric_Lebrun_Tod}), with the function $W$ and the one-form $A$ given by (\ref{W_equation1}) and (\ref{Equation_for_A2}) respectively, both of which are determined by $u$.

\subsection{An additional symmetry}

The $SU(\infty)$-Toda equation (\ref{Toda_3D}) is non-linear, and its general solution is not known. To progress further, we assume the existence of an additional symmetry. Tod \cite{Tod:2020ual} showed that without loss of generality, the corresponding Killing vector can be chosen to be $\partial_y$. We relabel the coordinate $y$ as $\phi$. The LeBrun--Tod ansatz becomes
\begin{equation}
	\label{metric_Lebrun_Tod2}
    \dif s^2=W^{-1}(\dif \tau+A)^2+We^u\dif\phi^2+W(\dif w^2+e^u\dif x^2),
\end{equation}
with $u$ satisfying the 2D $SU(\infty)$-Toda equation
\begin{equation}
	\label{Toda_equation_2D}
    (e^u)_{ww}+u_{xx}=0.
\end{equation}
The function $W$ satisfies the monopole equation (\ref{monopolo_equation}) without $\phi$- or $y$-dependence and is given by the same expression as (\ref{W_equation1}):
\begin{equation}
	\label{W_equation2}
	W=\frac{1}{k_1}\left(w-\frac{1}{2}w^2u_w\right).
\end{equation} 
With (\ref{W_equation2}) substituted in (\ref{Equation_for_A2}), the one-form $A$ now satisfies
\begin{equation}
	\dif A=\frac{1}{2k_1}\left[-w^2u_{xw}\dif w\wedge\dif \phi+e^u(w^2u_{ww}+w^2u_w^2-2wu_w+2)\dif x\wedge\dif \phi\right],
\end{equation}
which can be partially integrated to give
\begin{equation}
	\label{one_form_A0}
    A=\frac{1}{k_1}\left(B-\frac{1}{2}w^2u_x\right)\dif \phi, \quad \text{where } B_w=wu_x, \quad B_x=e^u(1-wu_w).
\end{equation}

\myheading{Weyl--Papapetrou coordinates}

With toric isometry, one-sided type-D Ricci-flat metrics can be cast in another canonical form in Weyl--Papapetrou coordinates. The toric isometry is generated by the Killing vectors $\partial_\tau$ and $\partial_\phi$. Consider the $2$-by-$2$ Gram-matrix $G_{ab}$, consisting of metric-tensor components $G_{ab}=g_{ab}$ with $a,b=\tau,\phi$. It is easy to see from (\ref{metric_Lebrun_Tod2}) that the determinant of $G_{ab}$ is
\begin{equation}
   \det (G_{ab})= e^u.
\end{equation}
The $SU(\infty)$-Toda equation satisfied by $u$ guarantees that a coordinate $\rho$ defined by
\begin{equation}
	\label{rho_coordinate}
    \rho \equiv e^{\frac{u}{2}}=\sqrt{\det(G_{ab})},
\end{equation}
is harmonic on the 2D base metric of (\ref{metric_Lebrun_Tod2})
\begin{equation}
	\label{metric_2D_base1}
    \dif s_{\text{2D}}^2=\dif w^2+e^u\dif x^2.
\end{equation}
From now on, the symbol $\rho$ shall be used as the coordinate defined in (\ref{rho_coordinate}), which bears no relation to the spin coefficient $\rho$. Denote the harmonic conjugate of $\rho $ as $z$ by
\begin{equation}
	\star_{\text{2D}}\dif\rho=\dif z,
\end{equation}
where $\star_{\text{2D}}$ is the Hodge star on the 2D base (\ref{metric_2D_base1}). The function $z$, which is used as another coordinate, is defined up to an additive constant. Flipping the orientation of the 2D base will flip the sign of $z$. Coordinates $(\rho,z)$ together with $(\tau,\phi)$ form the Weyl--Papapetrou coordinates in this context.

In the Weyl--Papapetrou coordinates,  the 2D base (\ref{metric_2D_base1}) is conformally mapped to the flat metric $\dif \rho^2+\dif z^2$. We use this result to cast the LeBrun--Tod ansatz (\ref{metric_Lebrun_Tod2}) in the Weyl--Papapetrou coordinates. We start with the 2D base metric (\ref{metric_2D_base1}) and transform it as follows: 
\begin{align}
	\dif s_{\text{2D}}^2&=\dif w^2+e^u\dif x^2=\dif w^2+\rho^2\dif x^2\nonumber\\
	&=(w_\rho^2 +\rho^2 x_\rho^2)\dif\rho^2+2(w_\rho w_z+\rho^2 x_\rho x_z)\dif\rho\,\dif z+(w_z^2+\rho^2 x_z^2)\dif z^2\nonumber\\
	&\propto \dif \rho^2+\dif z^2.
\end{align}
It follows that
\begin{equation}
	\begin{cases}
		w_\rho w_z+\rho^2 x_\rho x_z=0,\\
		w_\rho^2+\rho^2 x_\rho^2=w_z^2+\rho^2 x_z^2.
	\end{cases}
\end{equation}
Solving these equations for ($w_\rho,w_z$) or ($x_\rho,x_z$) yields
\begin{equation}
\label{w_and_x_relation}
    \begin{dcases}
    	w_\rho =\rho x_z\\
    	w_z=-\rho x_\rho
    \end{dcases}\quad\text{or}\quad
\begin{dcases}
	x_\rho =-\rho^{-1}w_z\\
	x_z=\rho ^{-1}w_\rho 
\end{dcases},
\end{equation}
up to a sign flip on one side of all these equations simultaneously. The same equations were obtained by Araneda and Lucietti \cite{Araneda:2025uqo} by inverting the Jacobian of the coordinate transformation from ($w,x$) to ($\rho,z$). The 2D base metric (\ref{metric_2D_base1}) in the Weyl--Papapetrou coordinates becomes
\begin{equation}
	\label{2D_quotient_in_Weyl}
    \dif s_{\text{2D}}^2=(w_\rho ^2 +w_z^2)(\dif \rho ^2 +\dif z^2)=\rho ^2 (x_\rho ^2 +x_z^2)(\dif \rho ^2 +\dif z^2).
\end{equation}

Now a key observation in our construction is that the equations for $x$ in (\ref{w_and_x_relation}) imply that $x$ satisfies the linear differential equation
\begin{equation}
\label{x_harmonicity}
    x_{\rho\rho}+\rho^{-1}x_\rho+x_{zz}=0,
\end{equation}
so it is an axisymmetric harmonic potential in an auxiliary 3D flat space with metric
\begin{equation}
	\label{3D_auxiliary_space1}
    \dif s_{\text{3D}}^2=\rho^2 \dif \phi ^2 +\dif\rho^2+\dif z^2.
\end{equation}
The equations for $w$ in (\ref{w_and_x_relation}) imply that $w$ is the ``harmonic conjugate'' of $x$. More explicitly, $w$ satisfies the linear differential equation
\begin{equation}
	\label{w_harmonicity}
    w_{\rho\rho}-\rho^{-1}w_\rho+w_{zz}=0,
\end{equation}
and it is an axisymmetric harmonic potential in an auxiliary 3D space with metric
\begin{equation}
	\label{3D_auxiliary_space2}
    \dif s_{\text{3D}}^2=\rho^{-2} \dif \phi ^2 +\dif\rho^2+\dif z^2,
\end{equation}
conjugate to (\ref{3D_auxiliary_space1}).

We record some useful calculation results that will be used later. Consider the Jacobian $\frac{\partial (w,x)}{\partial (\rho,z)}$ with components given by (\ref{w_and_x_relation}). We invert it to obtain
\begin{equation}
	\begin{dcases}
		\rho_w=\frac{w_\rho}{w_\rho^2+w_z^2}\\
		\rho_x=-\frac{\rho w_z}{w_\rho^2+w_z^2}
	\end{dcases},\qquad
	\begin{dcases}
	z_w=\frac{w_z}{w_\rho^2+w_z^2}\\
		z_x=\frac{\rho w_\rho}{w_\rho^2+w_z^2}	
	\end{dcases},
\end{equation}
where the right-hand sides of all these equations are expressed in terms of derivatives of $w$. Recall from (\ref{rho_coordinate}) that $u=2\ln\rho$. It then follows that
\begin{equation}
	\label{differential_of_u}
	\begin{dcases}
		u_w=\frac{2w_\rho}{\rho(w_\rho^2+w_z^2)},\\
		u_x=\frac{-2w_z}{w_\rho^2+w_z^2}.
	\end{dcases}
\end{equation}
By differentiating one more time, we can show that $u$ satisfies the 2D $SU(\infty)$-Toda equation (\ref{Toda_equation_2D}). The function $B$ in (\ref{one_form_A0}) now satisfies the following  differential equations
\begin{equation}
		\label{differential_of_B}
		\begin{dcases}
	B_\rho=B_w w_\rho+B_x x_\rho=-\rho w_z,\\ B_z=B_w w_z+B_x x_z=\rho w_\rho-2w.
	\end{dcases}
\end{equation}

\subsection{Toric one-sided type-D Ricci-flat metrics}

\label{sec_one_sided_type_D_metrics}

By substituting Eq.~(\ref{differential_of_u}) and Eq.~(\ref{differential_of_B}) to Eq.~(\ref{W_equation2}) and Eq.~(\ref{one_form_A0}), we can express the function $W$ and the one-form $A$ in terms of  ($\rho,z$). The results are then used to obtain the following form of toric one-sided type-D Ricci-flat metrics in the Weyl--Papapetrou coordinates $(\tau,\phi,\rho,z)$:
\begin{align}
	\label{metric_onesided_typeD}
    \dif s^2&=W^{-1}(\dif \tau+A)^2+W\rho^2\dif\phi^2+k_2^2WL(\dif \rho ^2 +\dif z^2),\\
    L&=w_\rho ^2 +w_z^2,\quad W=\frac{w}{k_1}\left(1-\frac{ww_\rho }{\rho L}\right),\quad 
    A=\frac{1}{k_1}\left(B+\frac{w^2w_z }{L}\right)\dif \phi.\label{one-form-A1}
\end{align}
The function $B$ in the one-form A satisfies
\begin{equation}
\label{B_equation}
B_\rho=-\rho w_z,\quad  B_z=\rho w_\rho-2w,
\end{equation}
and it is determined up to an additive constant, which can be absorbed by a redefinition of the coordinate $\tau$. The integrability of Eq.~(\ref{B_equation}) follows from Eq.~(\ref{w_harmonicity}) immediately.

The two constants $k_1$ and $k_2$ in the metric (\ref{metric_onesided_typeD}) can be set to 1 without loss of generality, but we prefer to keep them general, for later use to recover some known solutions more easily. In fact, in the metric (\ref{metric_onesided_typeD}), the constant factor $k_2^2$ is deliberately introduced in the 2D base metric conformal to $\dif \rho^2+\dif z^2$. This can be achieved as follows. First start with the metric (\ref{metric_onesided_typeD}) with $k_2=1$, and then rescale the coordinates $\tau$ and $\phi$ by doing the substitutions $(\tau\rightarrow k_2\tau,\phi\rightarrow k_2^{-1}\phi)$. Next we redefine parameter by doing the substitution $k_1\rightarrow k_1k_2^{-2}$. The metric is then brought to the same as (\ref{metric_onesided_typeD}) with the only change being the factor $k_2^2$ introduced. It is also clear that the fundamental two-form of the metric (\ref{metric_onesided_typeD}) is given by
\begin{equation}
	\label{Khaler_form}
	\omega=k_2\left[(\dif \tau+A)\wedge(w_\rho\dif\rho+w_z\dif z)+\rho W(-w_z \dif\rho+w_\rho\dif z)\wedge\dif\phi\right],
\end{equation}
and the Weyl scalar $\Psi_2$ is
\begin{equation}
	\label{Psi2_with_k2}
	\Psi_2=k_1k_2^{-2}w^{-3}.
\end{equation}

Thus, a one-sided type-D Ricci-flat solution in the Weyl--Papapetrou coordinates has the metric (\ref{metric_onesided_typeD}) and it is completely determined by the function $w$ that is axisymmetric and harmonic in the auxiliary 3D space (\ref{3D_auxiliary_space2}). Such a solution is also completely determined by the function $x$ even though it makes no explicit appearance in (\ref{metric_onesided_typeD}). The two functions $w$ and $x$ determine each other via (\ref{w_and_x_relation}). We shall refer to $w$ and $x$ as the generating potentials of the solution. The generating potential $x$ is defined in the auxiliary 3D flat space (\ref{3D_auxiliary_space1}) and has a more direct interpretation as describing Newtonian potentials produced by certain mass distributions.

The field equations for toric one-sided type-D Ricci-flat metrics have thus been linearised in the Weyl--Papapetrou coordinates, and they have been reduced to the pair of functions consisting of an axisymmetric harmonic potential $w$ in the auxiliary 3D space (\ref{3D_auxiliary_space2}) and its harmonic conjugate $x$. The harmonic potentials $(w,x)$ used here are nothing but the two coordinates that parametrise the 2D base in the LeBrun--Tod ansatz (\ref{metric_Lebrun_Tod2}), and $w$ is directly related to the Weyl scalar $\Psi_2$ (see (\ref{Psi2_with_k2})). It was Tod \cite{Tod:2020ual} who first performed this linearisation of field equations {and obtained} one-sided type-D Ricci-flat metrics with an additional symmetry. He introduced an axisymmetric harmonic potential $U$ in the auxiliary 3D flat space and applied a Ward transform
\begin{equation}
	\label{Ward_transform}
	{x=U_z,\quad w=\frac{1}{2}\rho U_\rho,}
\end{equation}
to linearise the field equations. {The integrability and harmonicity in the auxiliary 3D flat space of $U$ are guaranteed by Eq.~(\ref{w_and_x_relation}). Conversely, if $U$ is harmonic in the auxiliary 3D flat space, $x$ derived from the above Ward transform is harmonic and $w$ is the harmonic conjugate.} The harmonic potential $U$ introduced by Tod does not have a direct interpretation and its explicit form could be complicated even for simple solutions such as the Kerr solution \cite{Tod:2020ual}. {So, we propose to use the LeBrun--Tod coordinates $(w,x)$ as the generating potentials for toric one-sided type-D Ricci-flat metrics to simplify Tod's linearisation. The resulting form of metrics given by Eqs.~(\ref{metric_onesided_typeD})--(\ref{B_equation}) is simpler than the one originally given by Tod.}

Any pair of axisymmetric harmonic potentials $(w,x)$ as described above generates a one-sided type-D Ricci-flat metric. What ones generate one-sided type-D Ricci-flat multi-centre metrics? 

\section{Multi-centre solutions}

\label{sec_multi_soliton_solutions}

In the remainder of the paper, we shall refer to one-sided type-D Ricci-flat multi-centre metrics simply as multi-centre metrics or multi-centre solutions. We shall exclusively use axisymmetric harmonic potentials in this paper; therefore, we drop the term “axisymmetric,” and a harmonic potential will always mean an axisymmetric one.

\subsection{Harmonic potentials and gravitational solitons}

In this subsection, we focus on the generating potential $x$, which we recall is harmonic in the auxiliary 3D flat space (\ref{3D_auxiliary_space1}). We first introduce some functions, called gravitational solitons, on the half-plane defined by coordinates $(\rho,z)$ with $\rho\ge 0$ and use them to construct $x$ to generate the multi-centre solutions.

Gravitational solitons \cite{Belinski:2001ph} refer to those considered by Belinski and Zakharov \cite{Belinsky:1971nt} in their construction of vacuum solutions to Einstein's field equations using the inverse-scattering method (ISM). In the ISM, a dressing matrix is applied on a seed solution to generate a new solution. For solitonic solutions, the dressing matrix is representable as a rational function of the so-called spectral parameter with a finite number, say $n$, of simple poles $\mu_i$ ($1\le i\le n$). These simple poles are functions of the Weyl--Papapetrou coordinates $(\rho,z)$ with $\rho\ge 0$ and satisfy the following differential equations
\begin{equation}
\label{solition_PDE}
\mu_\rho=\frac{2\rho\mu}{\rho^2+\mu^2},\qquad\mu_z=-\frac{2\mu^2}{\rho^2+\mu^2}.
\end{equation}
The solutions to these equations are
\begin{equation}
\mu_i^\pm=\pm\sqrt{\rho^2+(z-z_i)^2}-(z-z_i),
\end{equation}
where $z_i$ is an integration constant. We call $\mu_i^+$ a gravitational soliton located at $z=z_i$ on the $z$-axis and $\mu_i^-$ an anti-soliton at the same position. Gravitational (anti-)solitons $\mu_i^\pm$ satisfy the following algebraic equations
\begin{equation}
\mu_i^+\mu_i^-=-\rho^2, \qquad\mu_i^++\mu_i^-=-2(z-z_i).
\end{equation}
They are the two roots of the quadratic equation $\mu^{2}+2(z-z_i)\mu^{2}-\rho^2=0$. Following \cite{Chen:2015iex}, we may, without loss of generality, restrict to solitons only. We denote $
\mu_i^+$ simply as $\mu_i$.
It is obviously non-negative on the $(\rho,z)$-plane with $\rho\ge 0$. An anti-soliton located at $z=z_i$ can be expressed in terms of a soliton at the same position as
\begin{equation}
\mu_i^-=-\frac{\rho^2}{\mu_i}.
\end{equation}

Gravitational solitons are closely related to harmonic potentials in the auxiliary 3D flat space (\ref{3D_auxiliary_space1}). In view of the equations (\ref{solition_PDE}), it follows that the logarithm of a soliton
\begin{equation}
x_i=\ln\mu_i, 
\end{equation} 
is a harmonic potential satisfying (\ref{x_harmonicity}). We say that $\ln \mu_i$ is the  harmonic potential associated with $\mu_i$. Near the $z$-axis, a soliton $\mu_i$ behaves as
\begin{align}
\mu_i\approx \begin{dcases}
2|z-z_i| &\text{as $\rho\rightarrow 0$\quad with fixed\quad }z< z_i,\\
\frac{\rho^2}{2|z-z_i|} &\text{as $\rho\rightarrow 0$ \quad with fixed\quad } z> z_i.
\end{dcases}
\end{align}
It is seen that the harmonic potential $\ln \mu_i$ diverges when we approach $\rho=0$ with fixed $z>z_i$. A closer look at this harmonic potential reveals that it is the Newtonian potential\footnote{We define the Newtonian potential produced by a point mass $m$ in the auxiliary 3D flat space at a distance $r$ from the mass as $\frac{m}{r}+\text{const.}$. {If the point mass is placed on the $z$-axis (as in our case) at $z=a$, then we have $r=\sqrt{\rho^2+(z-a)^2}$ in Weyl--Papapetrou coordinates.}} produced by a semi-infinite thin rod along the $z$-axis over $[z_i,\infty)$ with unit linear mass density in the auxiliary 3D flat space. With the presence of multiple solitons, we can add or subtract their associated harmonic potentials and interpret the resulting potential as being produced by a thin rod along the $z$-axis with piecewise uniform linear mass density. A few more examples are listed here:
\begin{equation}
\label{rod_source_interpretation}
\begin{aligned}
\ln\rho^2&:\text{Infinite rod with unit linear mass density over } (-\infty,\infty);\\
\ln\mu_i&: \text{Semi-infinite rod  with unit linear mass density over } [z_i,\infty);\\
\ln\frac{\rho^2}{\mu_i}&:\text{Semi-infinite rod  with unit linear mass density over } (-\infty,z_i];\\
\ln \frac{\mu_i}{\mu_j}&:\text{Finite rod  with unit linear mass density over } [z_i,z_j].
\end{aligned}
\end{equation}

Gravitational solitons and their associated harmonic potentials have been used in the construction of static and axisymmetric solutions in Lorentzian regime, known as Weyl's class of metrics. Einstein's field equations in this case reduce to a single harmonic potential in the same auxiliary 3D flat space (\ref{3D_auxiliary_space1}), and the logarithms of metric-tensor components $\ln |g_{tt}|$ and $\ln g_{\phi\phi}$ are linear combinations of harmonic potentials associated with gravitational solitons as discussed above, see, e.g., \cite{Emparan:2001wk}. Here, we have assumed that the time translational symmetry and the axial symmetry are generated by $\partial_t$ and $\partial_\phi$ respectively. The Schwarzschild solution and the static C-metric belong to Weyl's class of metrics.  Now toric one-sided type-D Ricci-flat metrics are also reduced to a single harmonic potential $x$ in the auxiliary 3D flat space, even though the algorithm to rebuild the metrics from the harmonic potential in this case is very different from that for Weyl's class of metrics. Nevertheless, we shall attempt to use harmonic potentials (\ref{rod_source_interpretation}) and their linear combinations to construct the potential $x$ to generate one-sided type-D Ricci-flat metrics.

\subsection{Multi-centre solutions}

\myheading{The generating potential $x$}

To generate the multi-centre solutions, we propose the following form of the generating potential $x$:
\begin{equation}
	\label{x_general1}
	x=\alpha_0\ln \rho+\sum_{i=1}^{n} \alpha_i\ln \mu_i,
\end{equation}
where all $\alpha$'s are constants. We observe that for the Kerr--NUT solution, the generating potential $x$ indeed has this form with $n=2$,\footnote{The generating potential $x$ for the Kerr--NUT solution can be found by identifying $L=\rho^2(x_\rho^2+x_z^2)$ in (\ref{one-form-A1}) with $\frac{F}{R_{11}R_{22}}$ in (5.13) in \cite{Chen:2015iex} and solving for $x$.} and moreover, it satisfies $\alpha_0=-(\alpha_1+\alpha_2)$. We thus further propose that 
\begin{equation}
	\alpha_0=-\sum_{i=1}^{n}\alpha_i.
\end{equation}
This choice of $\alpha_0$ is also necessary to ensure that the second factor of $W$ in (\ref{one-form-A1}) does not have zeros on and away from the $z$-axis (for $\rho\ge 0$), failure of which will result in curvature singularities at the zeros. So our choice of the generating potential $x$ for the multi-centre solutions is
\begin{equation}
	\label{x_of_multisoliton}
	x=\sum_{i=1}^{n} \alpha_i\ln \fracb{\mu_i}{\rho}.
\end{equation}
We remark that the same choice of $\alpha_0$ was made in earlier works in \cite{Biquard:2021gwj,Araneda:2025uqo}, to ensure that the resulting solutions have either ALF or ALE asymptotic structure.

\myheading{The generating potential $w$}

With the generating potential $x$ given by (\ref{x_of_multisoliton}), we calculate $w_\rho$ and $w_z$ using (\ref{w_and_x_relation}) to get
\begin{align}
w_\rho&=\rho 	x_z=-\sum_{i=1}^{n}\frac{\alpha_i\rho}{\sqrt{\rho^2+(z-z_i)^2}}=-\sum_{i=1}^{n}\frac{2\alpha_i\rho\mu_i}{\rho^2+\mu_i^2},\label{w_rho}\\
w_z&=-\rho x_\rho=-\sum_{i=1}^{n}\fracb{\alpha_i(z-z_i)}{\sqrt{\rho^2+(z-z_i)^2}}=-\sum_{i=1}^{n}\frac{\alpha_i(\rho^2-\mu_i^2)}{\rho^2+\mu_i^2}.\label{w_z}
\end{align}
The generating potential $w$ can be integrated as
\begin{equation}
	\label{w_of_multisoliton}
	w=-\left[\beta+\sum_{i=1}^{n}\alpha_i\sqrt{\rho^2+(z-z_i)^2}\right]=-\left[\beta+\sum_{i=1}^{n}\frac{\alpha_i (\rho^2+\mu_i^2)}{2\mu_i}\right],
\end{equation}
where $\beta$ is an integration constant, which turns out to be important to determine the asymptotic structure of the solution. The function $B$ in (\ref{B_equation}) can be integrated as
\begin{equation}
	\label{B_function}
	B=2\beta z+\sum_{i=1}^{n}\alpha_i(z-z_i)\sqrt{\rho^2+(z-z_i)^2}=2\beta z+\sum_{i=1}^{n}\frac{\alpha_i}{4}\left(\frac{\rho^4}{\mu_i^2}-\mu_i^2\right),	
\end{equation}
where an integration constant has been omitted. In (\ref{B_function}), the function $z$ following the second equal sign can be replaced by the expression $\frac{\rho^2-\mu^2}{2\mu}$ for a soliton $\mu$ with an arbitrary position.

\myheading{Multi-centre solutions}

It will be shown in Sec.~\ref{rod_structre} that, with the generating potentials $(w,x)$ given by (\ref{w_of_multisoliton}) and (\ref{x_of_multisoliton}), the function $W$ in the metric (\ref{metric_onesided_typeD}) has a pole at each $z=z_i$, the position of the soliton $\mu_i$, along the $z$-axis. So $z=z_i$ represents a centre of the solution. Since there are $n$ solitons present in the harmonic potentials (\ref{w_of_multisoliton}) and (\ref{x_of_multisoliton}), there are $n$ centres in the solution. One-sided type-D Ricci-flat multi-centre solutions are then given by Eqs.~(\ref{metric_onesided_typeD})--(\ref{B_equation}), with various expressions in these equations specified by Eqs.~(\ref{w_rho})--(\ref{B_function}). They are the one-sided type-D counterpart of the Gibbons--Hawking multi-centre metrics {with all centres collinearly aligned.} We note that, for our multi-centre solutions, all metric-tensor components can be written as rational functions of the various solitons $\mu_i$ and the expression $\rho^2$, a property that was also observed for solutions constructed in \cite{Chen:2015iex}.

When $n=1$, we have the 1-centre solution; when $n=2$, we have the 2-centre solution, etc.. The parameter $z_i$ is the position of the $i$-th centre, and $\alpha_i$ will be called the charge carried by the $i$-th centre. When the charge $\alpha_i=0$, the centre $z=z_i$ disappears from the solution: $z=z_i$ is not a pole of $W$ anymore, and the corresponding soliton $\mu_i$ drops out from the solution---we have one centre less. The fundamental two-form (\ref{Khaler_form}) specifies the orientation of the $n$-centre solution.  By construction, the Petrov type of the $n$-centre solution is D$^+$X$^-$ for some X $\in\{\text{I, D, O}\}$. The Petrov type X of the anti-self-dual Weyl tensor depends on the number of centres $n$ (and also other parameters in the solution) and has to be computed case by case.  

We remark that the multi-centre solutions described above appeared in equivalent forms in earlier works by Biquard and Gauduchon \cite{Biquard:2021gwj,Biquard:2023rbr}. {If we perform the Ward transform (\ref{Ward_transform}) on the harmonic potential $x$  given by (\ref{x_of_multisoliton}) (with $w$ calculated as Eq.~(\ref{w_of_multisoliton})), we can integrate to find Tod's generating potential $U$ for the multi-centre solutions as
\begin{equation}
	U=\sum_{i=1}^{n}\alpha_i\left[(z-z_i) \ln \frac{\mu_i}{\rho}+\frac{\rho^2+\mu_i^2}{2\mu_i}\right]+C,
\end{equation}
where $C$ is an arbitrary integration constant. A direct interpretation of this potential is not yet known. Biquard and Gauduchon considered Tod's ansatz with the above generating potential in an equivalent form, and, by imposing regularity conditions, extracted and classified all toric one-sided type-D gravitational instantons. Our focus in this paper, on the other hand, is on the multi-centre metrics themselves. In particular, we allow mild singularities such as conical and orbifold ones, so the multi-centre metrics we study do not necessarily represent gravitational instantons. Nevertheless, the classification of these metrics is still a very interesting problem.}

\myheading{Parameter counting}

Here, we count the number of independent parameters of the $n$-centre solution given by Eqs.~(\ref{metric_onesided_typeD})--(\ref{B_equation}) and Eqs.~(\ref{w_rho})--(\ref{B_function}). We note that the $n$-centre solution is invariant under a few operations, called symmetries of the solution.

\begin{enumerate}
	\item Translational symmetry along the $z$-axis 
	
	The $n$-centre solution is invariant under the following simultaneous substitutions
	\begin{equation}
		\label{symmetry_translational}
		z_i\rightarrow z_i+c_1,\quad z\rightarrow z+c_1,\quad \tau\rightarrow \tau-\frac{2\beta c_1}{k_1}\phi,
	\end{equation}
	for an arbitrary constant $c_1$. The redefinition of $\tau$ is necessary since a constant $2\beta c_1$ will be generated in (\ref{B_function}) under the above substitutions.

	\item Parameter-rescaling symmetry
	
	The $n$-centre solution is invariant under the following simultaneous substitutions
	\begin{equation}
		\label{symmetry_scaling}
		\alpha_i\rightarrow c_2
		\alpha_i,\quad \beta\rightarrow c_2\beta,\quad k_1\rightarrow c_2k_1, \quad k_2\rightarrow c_2^{-1}k_2,
	\end{equation}
	for any non-zero constant $c_2$.
	
	\item Coordinate-rescaling symmetry on ($\tau,\phi$)
	
	The $n$-centre solution is invariant under the following simultaneous substitutions
	\begin{equation}
		\tau\rightarrow c_3\tau,\quad \phi\rightarrow c_3^{-1}\phi,\quad k_1\rightarrow c_3^{-2}k_1,\quad k_2\rightarrow c_3^{-1}k_2,
	\end{equation}
	for any constant $c_3>0$. Flipping the sign of $c_3$ does not change the metric. We will only consider positive $c_3$.
	
	\item Coordinate-rescaling symmetry on $(\rho,z)$
	
	The $n$-centre solution is invariant is under the following simultaneous substitutions
	\begin{equation}
		(\rho,z,z_i)\rightarrow c_4(\rho,z,z_i),\quad \phi\rightarrow c_4^{-1}\phi,\quad \beta\rightarrow c_4\beta,\quad k_1\rightarrow c_4k_1\quad k_2\rightarrow c_4^{-1}k_2,
	\end{equation}
	for any constant $c_4>0$. The positivity of $c_4$ is necessary for $\rho$ to remain positive as preferred.
\end{enumerate}

The translational symmetry implies that among the parameters $z_i$, only $n-1$ of them are independent. The parameter-rescaling symmetry implies that among the parameters $(\beta,\alpha_i)$, only $n$ of them are independent. It also implies that $k_1$ can be chosen to be positive without loss of generality, which is what we will do in the following. We will also choose $k_2>0$, which is always possible since only $k_2^2$ appear in the metric. The coordinate-rescaling symmetries 3 and 4 then imply that we can set $k_1$ and $k_2$ to any positive values without loss of generality. The $n$-centre solution thus has a total of $2n-1$ independent parameters: $n-1$ parameters among $z_i$ and another $n$ parameters among $(\beta,\alpha_i)$. The 1-, 2- and 3-centre solutions have 1, 3, and 5 independent parameters respectively. This counting of independent parameters agrees with that found in explicit metrics constructed in the next section.

Since all centres are treated on an equal footing, we can assume without loss of generality the ordering of centre positions
\begin{equation}
	z_1<z_2<...<z_{n-1}<z_n.
\end{equation}
It will be shown in the following subsection that, with $k_1$ chosen as positive, it is necessary that $\beta\ge 0$ and $\alpha_i\ge 0$ for all $1\le i\le n$ 
with not all the $\alpha$'s being zero at the same time. Different normalisations can be imposed on the parameters $(\beta,\alpha_i)$ using the parameter-rescaling symmetry. If we fix the parameter $\beta=1$, the multi-centre solutions are completely determined by the positions $z_i$ and charges $\alpha_i$ of the centres, in much the same way that the multi-centre Gibbons--Hawking metrics are completely determined by the positions and charges of centres. Another possibility is to impose
$\sum\limits_{i=1}^n\alpha_i=1$ and leave $\beta$ free. In this paper, different normalisation conditions may be imposed for different solutions (e.g., solutions with different $n$ and different asymptotic structures), on a case-by-case basis, to recover known forms of metrics that are used in the literature.

\subsection{Rod structure}

\label{rod_structre}

For clarity of presentation, we first introduce ``flip-sign-from-$k$'' summation notation $\ksum$. With a given positive integer $n$, we define
\begin{equation}
	\bigksum h_i = \sum\limits_{1\le i\le k-1}h_i-\sum\limits_{ k\le i\le n}h_i,
\end{equation}
which carries an integer parameter $k$ satisfying $1\le k\le n+1$ and the sum is done over the index $i$ from $1$ to $n$ with a minus sign inserted for terms with $i\ge k$. We will omit summation limits when a (usual) summation runs from $i=1$ to $i=n$, i.e.,
\begin{equation}
	\sum h_i\equiv\sum\limits_{i=1}^n h_i.
\end{equation}
The following summation identities can be easily proved:
\begin{align}
	&\bigksum h_i=\begin{dcases}
		-\sum h_i &\text{ for $k=1$}\\
		\sum h_i &\text{ for $k=n+1$}
	\end{dcases},\\	
	&\bigksum (h_i+g_i)=\bigksum h_i+\bigksum g_i,\label{k-sum-linearity}\\
	&	\kplusonesum h_i=\bigksum h_i+2h_k.\label{k-sum-difference}
\end{align}

\myheading{Rod source of the generating potential $x$}

Consider the generating potential $x$ given by (\ref{x_of_multisoliton}) of the $n$-centre solution. We shall refer to the mass distribution in the auxiliary 3D flat space that is responsible to produce a Newtonian potential equal to $x$ as the rod source of the solution. We first note that in (\ref{x_of_multisoliton}) we can write the expression $\ln \frac{\mu_i}{\rho}$ as
\begin{equation}
	\ln \frac{\mu_i}{\rho}=\frac{1}{2}\ln \mu_i-\frac{1}{2}\ln \frac{\rho^2}{\mu_i}.
\end{equation}
According to our interpretation (\ref{rod_source_interpretation}), this is the potential produced by an infinite thin rod along the $z$-axis that has linear mass density $\lambda=-\frac{1}{2}$ on $(-\infty,z_i]$ and $\lambda=\frac{1}{2}$ on $[z_i,\infty)$, with a jump of one unit in $\lambda$ when $z=z_i$ is crossed. The generating potential $x$ given in (\ref{x_of_multisoliton}) is seen as the potential produced by a thin rod along the $z$-axis with piecewise uniform linear mass density, whose value jumps by an amount of $\alpha_i$ when $z=z_i$ is crossed. So, there are in total $n+1$ rod segments in the rod source, each with a uniform linear mass density. The $k$-th rod, or rod $k$, with $z\in [z_{k-1},z_k]$ has linear mass density
\begin{equation}
	\label{liear_mass_dentisity_k}
	\lambda_k=\frac{1}{2}\ksum\alpha_i.
\end{equation}
The first rod has $z\in(z_0\equiv-\infty,z_1]$ and $\lambda_1=-\frac{1}{2}\sum\alpha_i$; the last rod has $z\in[z_n,z_{n+1}\equiv\infty)$ and $\lambda_{n+1}=\frac{1}{2}\sum\alpha_i$.

A few remarks are in order. First, here we allow rod segments with negative linear mass density. Second, the rod source is the mass distribution responsible to produce the generating potential $x$, not any metric-tensor components themselves. This is to be distinguished from the usage of rod source in \cite{Emparan:2001wk} and \cite{Harmark:2004rm} in the study of static/stationary and axisymmetric black hole solutions. Third, the linear mass density is ill-defined and can be assigned with any finite value at a rod junction. The potential generated will not be affected since the rod junctions are isolated points with zero measure.  Finally, all centre positions $z_i$ and charges $\alpha_i$ are encoded in the rod source by the positions of rod junctions and linear mass densities of various rod segments. As discussed earlier, if the normalisation condition $\beta=1$ is chosen, these quantities are independent parameters and they completely determine the $n$-centre solution.  

\myheading{Poles of $W$ along the rod axis}

Now we shift our attention from the generating potential $x$ to metric-tensor components, and focus on the near-rod behaviour of $g_{\tau\tau}^{-1}=W$. We first note the following leading-order behaviours of various functions (\ref{w_rho})--(\ref{B_function}) in the metric (\ref{metric_onesided_typeD}) near the interior of rod $k$ for $(\rho\rightarrow 0,z_{k-1}<z<z_{k})$:
\begin{equation}
	\begin{aligned}
		w_\rho&\approx -\ksum\frac{\alpha_i\rho}{z-z_i},\quad
		w_z\approx-\ksum\alpha_i,\quad w\approx -\beta -\ksum \alpha_i(z-z_i),\\
		L&\approx\left(\ksum \alpha_i\right)^2,\quad
		B\approx 2\beta z+\ksum \alpha_i(z-z_i)^2.
	\end{aligned}\label{rod_functions_behavior}
\end{equation}
In the interior of rod $k$ for $(\rho= 0,z_{k-1}<z<z_{k})$, the function $W=g_{\tau\tau}^{-1}$ behaves as
\begin{equation}
	\label{W_onrod}
	W|_{\text{rod k}}= \frac{1}{k_1}\left(\beta+\ksum\alpha_i(z-z_i)\right)\left[\frac{(\beta+\ksum \alpha_i(z-z_i))}{(\ksum\alpha_i)^2}\ksum\frac{\alpha_i}{z-z_i}-1\right],
\end{equation}
{which is a rational function in terms of $z$.} We see that the end points of a rod correspond to poles of $W$ (as a function of $z$). Using (\ref{W_onrod}), we can find the behaviour of $W$ around the pole $z=z_k$ as
\begin{align}
	W\approx\begin{dcases}
		\frac{1}{k_1}\left(\frac{\beta+\bigksum\alpha_i(z_k-z_i)}{\bigksum \alpha_i}\right)^2\left(-\frac{\alpha_k}{z-z_k}\right) \quad&\text{for} \quad z\rightarrow z_k^-;\\
		\frac{1}{k_1}\left(\frac{\beta+\kplusonesum\alpha_i(z_k-z_i)}{\kplusonesum \alpha_i}\right)^2\left(\frac{\alpha_k}{z-z_k}\right) \quad&\text{for} \quad z\rightarrow z_k^+.
	\end{dcases}
\end{align}
It is thus clear that, with $k_1$ chosen as positive, all the $\alpha$'s should be non-negative to {ensure $W>0$ around all the poles}:
\begin{equation}
	\label{range_alpha}
	\alpha_i\ge 0 \quad (i=1,2,...,n).
\end{equation} 
Now examine the behaviour (\ref{W_onrod}) of $W$ on the last rod (rod $n+1$) with $z\in[z_n,\infty)$ and $z\rightarrow \infty$, and we find that
\begin{equation}
	W\approx \frac{\beta}{k_1} .
\end{equation}
Hence, {for $W$ to be positive along this rod}, it is also required that
\begin{equation}
	\label{range_beta}
	\beta\ge 0.
\end{equation}
{With the above choice of ranges of parameters for $\alpha_i$ and $\beta$, the function $W$ can be seen to be smooth and strictly positive in the interior of each rod, including the two rods that extend to infinity.}

Now we proceed to show that the ranges (\ref{range_alpha}) and (\ref{range_beta}) are sufficient to ensure that $W$ is positive away from the rod axis for $\rho>0$. We first define a distance function $r_i$ for each $1\le i\le n$ as
\begin{equation}
	\label{distance_function}
	r_i\equiv\sqrt{\rho^2+(z-z_i)^2},
\end{equation}
which measures distance  in the auxiliary 3D flat space from the point under consideration to the $i$-th centre. The generating potential $w$ given by (\ref{w_of_multisoliton}) can be written as
\begin{equation}
	w=-\beta-\sum\alpha_ir_i,
\end{equation}
which is obviously negative. In view of the second equation in (\ref{one-form-A1}), to show the positivity of $W$ for $\rho>0$, it is sufficient to show that $1-\frac{w w_\rho}{\rho L}<0$. The proof is presented below:
\begin{align}
	&1-\frac{w w_\rho}{\rho L}=1-\frac{(\beta+\sum\alpha_ir_i)\sum \frac{\alpha_i}{r_i}}{\left(\sum\frac{\alpha_i(z-z_i)}{r_i}\right)^2+\left(\sum\frac{\alpha_i \rho}{r_i}\right)^2}<0\nonumber\\
	\Leftarrow\quad&\left(\sum\frac{\alpha_i(z-z_i)}{r_i}\right)^2+\left(\sum\frac{\alpha_i\rho}{r_i}\right)^2<\left(\beta+\sum\alpha_ir_i\right)\sum \frac{\alpha_i}{r_i}\nonumber\\
	\Leftarrow\quad&\left(\sum\alpha_i\cos\theta_i\right)^2+\left(\sum\alpha_i\sin\theta_i\right)^2<\sum\alpha_i r_i\sum\frac{\alpha_i}{r_i}\nonumber\\
	\Leftarrow\quad &\sum\alpha_i^2+2\sum_{i<j}\alpha_i\alpha_j\cos(\theta_i-\theta_j)<\sum\alpha_i^2+\sum_{i<j}\alpha_i\alpha_j\left(\frac{r_i}{r_j}+\frac{r_j}{r_i}\right)\nonumber
\end{align}
In the third line we have introduced $\cos\theta_i=\frac{z-z_i}{r_i}$ and $\sin\theta_i=\frac{\rho}{r_i}$ and omitted $\beta$ in the estimate. The last line is true since all the $r$'s are positive and distinct.

It then follows that $g_{\tau\tau}>0$ away from the rod axis. We also have $\det g_{ab}=\rho^2>0$ for $a,b=\tau,\,\phi$ and $g_{\rho\rho}=g_{zz}=k_2^2WL>0$ away from the rod axis. The metric thus has Euclidean signature everywhere away from the rod axis in the half-plane with $(\rho>0,-\infty<z<\infty)$.

\myheading{Rod directions}

Toric Ricci-flat metrics can be analysed using the so-called rod-structure formalism. Several versions of it have been developed in the literature. In this paper we follow \cite{Chen:2010zu}. Consider a toric Ricci-flat metric written in the Weyl--Papapetrou coordinates in the form
\begin{equation}
	\label{metric_toric}
	\dif s^2=g_{\tau\tau}(\dif \tau+A_\phi\dif\phi)^2+g_{\tau\tau}^{-1}\rho^2\dif\phi^2+f(\dif\rho^2+\dif z^2).
\end{equation}
The $2\times 2$ Gram-matrix $G_{ab}=g_{ab}$ for $a,b=\tau,\phi$ has determinant $\det G_{ab}=\rho^2$. So it has at least one-dimensional kernel along the $z$-axis with $\rho=0$. At isolated points $z=z_i$ along the $z$-axis, the Gram-matrix $G_{ab}$ vanishes and has a two-dimensional kernel. These points on the $z$-axis are called turning points of the solution. Segments on the $z$-axis between two adjacent turning points are called rods. If there are $n$ turning points in the solution, there will be $n+1$ rods: $(-\infty,z_1],[z_1,z_2],...,[z_n,\infty)$, labeled as rod 1, rod 2, etc.. In the interior of each rod, the Gram-matrix $G_{ab}$ has a kernel that is exactly one-dimensional. For the interior of rod $k$, a vector $v^0_k$ in this kernel is called an unnormalised direction of that rod. For the metric given in (\ref{metric_toric}), an unnormalised rod direction $v^0_k$ and the corresponding surface gravity $\kappa_k$ for rod $k$ is given by (see \cite{Chen:2015iex})
\begin{equation}
	\label{rod_direction_unnormalised}
	v^0_k=(-A_\phi,1)|_{\text{rod k}}, \quad \kappa_k^{-1}=\sqrt{fg_{\tau\tau}}|_{\text{rod k}},
\end{equation}
where evaluation of the above expressions is taken in the interior of rod $k$ with $(\rho=0, z_{k-1}<z<z_k)$. The surface gravity $\kappa_k$ can be used to define the normalised rod direction $v_k$ (up to a sign) for rod $k$ as
\begin{equation}
	\label{rod_direction_normalised}
	v_k=\kappa_k^{-1}v^0_k.
\end{equation}
Both the unnormalised and normalised rod directions (\ref{rod_direction_unnormalised}) and (\ref{rod_direction_normalised}) are independent of $z$ in the interior of each rod. The metric under study is free of conical singularity along rod $k$ if $v_k$ generates flow of period of $2\pi$. The specification of all turning points $z_i$ $(1\le i\le n)$ and normalised rod directions $v_k$ $(1\le k\le n+1)$ forms the rod structure of a solution. 

The rod structure of the $n$-centre solution given by Eqs.~(\ref{metric_onesided_typeD})--(\ref{B_equation}) and Eqs.~(\ref{w_rho})--(\ref{B_function}) can be calculated directly. There are $n$ turning points, corresponding to the $n$ centre positions, i.e., $z=z_i$ $(1\le i\le n)$, which divide the $z$-axis in $n+1$ rods. The metric (\ref{metric_onesided_typeD}), when cast in the form of (\ref{metric_toric}), has components
\begin{equation}
	\label{metric_tensor_components}
	g_{\tau\tau}=W^{-1}, \quad f=k_2^2WL,
\end{equation}
and $A_\phi$ is determined from the one-form $A=A_\phi\dif\phi$ given in (\ref{one-form-A1}). For rod $k$, using the expressions in (\ref{rod_functions_behavior}) and formulas (\ref{rod_direction_unnormalised})--(\ref{metric_tensor_components}), direct calculation yields
\begin{align}
	\kappa_k^{-1}&=k_2\sqrt{L}|_{\text{rod k}}=k_2\ksum\alpha_i,\\
	A_\phi |_{\text{rod k}}&=-\frac{1}{k_1}\left[\frac{\left(-\beta+\ksum\alpha_iz_i\right)^2}{\ksum\alpha_i}-\ksum\alpha_i z_i^2\right].
\end{align}
The flip-sign-from-$k$ summation identity (\ref{k-sum-linearity}) has been used to simplify calculations. The normalised direction $v_k$ of rod $k$ is thus
\begin{equation}
	\label{rod_directions}
	v_k\equiv(a_k,b_k)=\left(\frac{k_2}{k_1}\left[\left(-\beta+\ksum\alpha_i z_i\right)^2-\ksum\alpha_i\ksum\alpha_i z_i^2\right],k_2\ksum\alpha_i\right).
\end{equation}
Rod directions from this point onwards refer to normalised ones.

Some properties of the rod directions are worth mentioning. We recall that for the $n$-centre solution, each rod has a uniform linear mass density. For rod $k$, we find that the second component $v_k[2]=b_k$ of its normalised direction $v_k$ is proportional to its linear mass density $\lambda_k$ (cf. (\ref{liear_mass_dentisity_k}))
\begin{equation}
	\label{b_k}
	b_k=2k_2\lambda_k.
\end{equation}
Using flip-sign-from-$k$ identity (\ref{k-sum-difference}), we can deduce that
\begin{equation}
	\label{jump_in_b}
b_{k+1}-b_k=2k_2\alpha_k>0.
\end{equation} 
It then follows that $b_k$ ($1\le k\le n+1$) for different rods satisfy
\begin{align}
	b_1<b_2<...<b_{n}<b_{n+1}.
\end{align}
We note that $b_1=-k_2\sum \alpha_i<0$ and $b_{n+1}=k_2\sum\alpha_i>0$. Define $2\times 2$ matrices $M_k$ $(1\le k\le n)$ formed by the ordered pair $(v_k,v_{k+1})$ of adjacent rod directions as
\begin{equation}
	M_k\equiv (v_k,v_{k+1})=\left(\begin{matrix}
		a_k&a_{k+1}\\b_k&b_{k+1}
	\end{matrix}\right).
\end{equation} 
Using flip-sign-from-$k$ summation identities (\ref{k-sum-linearity}) and (\ref{k-sum-difference}), the determinant of this matrix can be calculated as
\begin{equation}
	\label{M_K}
	\det M_k=\frac{2k_2^2\alpha_k}{k_1}\left(\beta+\ksum\alpha_i(z_k-z_i)\right)^2,
\end{equation}
which is seen to be strictly positive.

\subsection{Asymptotic structure}

In this subsection, we show explicitly that the $n$-centre solution is ALF. The infinity is located at $\sqrt{\rho^2+z^2}\rightarrow \infty$. If we introduce coordinates $(r,\theta)$ by
\begin{equation}
	\label{CT_asymptotic_structure}
	\rho=r\sin\theta,\qquad z=r\cos\theta,
\end{equation}
then infinity is located at $r\rightarrow \infty$.

Near infinity at large $r$, the various functions (\ref{w_rho})--(\ref{B_function}) that appear in the metric (\ref{metric_onesided_typeD}) behave as follows:
\begin{equation}
	\begin{aligned}
		w_\rho&\approx-\sin\theta\sum\alpha_i-\frac{\sin\theta\cos\theta}{r}\sum\alpha_iz_i-\frac{\sin\theta(3\cos^2\theta-1)}{2r^2}\sum\alpha_iz_i^2,\\
		w_z&\approx-\cos\theta\sum\alpha_i+\frac{\sin^2\theta}{r}\sum\alpha_iz_i+\frac{3\sin^2\theta\cos\theta}{2r^2}\sum\alpha_i z_i^2,\\
		w&\approx-\left(r\sum\alpha_i+\beta-\cos\theta\sum\alpha_iz_i+\fracb{\sin^2\theta}{2r}\sum\alpha_iz_i^2\right),\\
		L&\approx\left(\sum\alpha_i\right)^2-\fracb{\sin^2\theta}{r^2}\left[\sum\alpha_i\sum\alpha_i z_i^2-\left(\sum\alpha_iz_i\right)^2\right],\\
		B&\approx r^2\cos\theta\sum\alpha_i+2\beta r\cos\theta-r(1+\cos^2\theta)\sum\alpha_iz_i\\
		&\,\,\,\,+\cos\theta\left(1+\frac{\sin^2\theta}{2}\right)\sum\alpha_iz_i^2,
	\end{aligned}
\end{equation}
where higher-order terms have been omitted. It follows that at large $r$
\begin{align}
	g_{\tau\tau}&\approx\frac{k_1}{\beta},\qquad
	g_{\rho\rho}\approx\frac{\beta}{k_1}\left(k_2\sum\alpha_i\right)^2,\\
	A&\approx\frac{1}{k_1\sum\alpha_1}\left[\left(\sum\alpha_i\sum\alpha_i z_i^2-\left(\sum\alpha_i z_i\right)^2-\beta^2\right)\cos\theta+2\beta\sum\alpha_i z_i\right]\dif\phi.\label{one_form_infinity}
\end{align}
The coordinate transformation (\ref{CT_asymptotic_structure}) gives $\dif \rho^2+\dif z^2=\dif r^2 +r^2 \dif \theta^2$. So if we fix the constants $k_1$ and $k_2$ as
\begin{equation}
	\label{k12}
	k_1=\beta, \qquad k_2=\left(\sum\alpha_i\right)^{-1},
\end{equation}
respectively, we will have $g_{\tau\tau}\approx1$ and $g_{\rho\rho}=g_{rr}\approx1$ at infinity. The one-form $A$ in (\ref{one_form_infinity}) has the form
\begin{equation}
	A\approx\left(2n_\text{nut}\cos\theta+\frac{2\sum \alpha_iz_i}{\sum\alpha_i}\right)\dif\phi,
\end{equation}
with $2n_\text{nut}$ given by
\begin{align}
	\label{NUT-charge}
	2n_\text{nut}=\frac{1}{\sum\alpha_i}\left(\frac{\sum\alpha_i\sum\alpha_i z_i^2-\left(\sum\alpha_i z_i\right)^2}{\beta}-\beta\right).
\end{align}
With the constant term in the one-form $A$ dropped, at large $r$, the metric (\ref{metric_onesided_typeD}) approaches
\begin{align}
	\label{metric_ALF}
	ds^2\approx \left(\dif\tau+2n_\text{nut}\cos\theta\dif\phi\right)^2+\dif r^2+r^2\left(\dif \theta^2+\sin^2\theta\dif\phi^2\right).
\end{align}
We recognise that $n_\text{nut}$ is the asymptotic NUT charge. The Killing vector $\partial_\tau$ has unit norm at infinity. Further fall-off behaviour of the metric-tensor components near infinity can be determined, but we will not pursue here. The $n$-centre solution is seen to be ALF.

Note that $\sum\alpha_i\sum\alpha_i z_i^2-\left(\sum\alpha_i z_i\right)^2> 0$ for $n\ge 2$, so the first term in (\ref{NUT-charge}) is positive. The NUT charge $n_{\text{nut}}$ given by (\ref{NUT-charge}) is a strictly decreasing function of $\beta$ when all centre positions $z_i$ and charges $\alpha_i$ are fixed. We also note that
\begin{equation}
	\label{NUT_charge_infinity}
	\begin{cases}
		n_{\text{nut}}\rightarrow +\infty\text{ as } \beta\rightarrow 0,\\
		n_{\text{nut}}\rightarrow -\infty \text{ as }\beta\rightarrow +\infty.
	\end{cases}
\end{equation}
So any $n_{\text{nut}}\in \mathbb{R}$ can be achieved by choosing $\beta$ appropriately. The $n$-centre solution is asymptotically flat (AF) when
\begin{equation}
	n_{\text{nut}}=0\quad \Leftrightarrow\quad \beta=\sqrt{\sum\alpha_i\sum\alpha_i z_i^2-\left(\sum\alpha_i z_i\right)^2}.
\end{equation}

An exceptional case happens when $n=1$, the first term in the NUT charge $n_{\text{nut}}$ given by (\ref{NUT-charge}) vanishes so
\begin{equation}
	2n_{\text{nut}}=-\frac{\beta}{\alpha_1}\in (-\infty,0).
\end{equation}
The NUT charge is always negative for the 1-centre solution. The AF limit in this case is the trivial limit with $\beta=0$.

\subsection{Two ALE limits}

The $n$-centre solution is ALF, with a compact dimension at infinity whose size is determined by the NUT charge (\ref{NUT-charge}). In the limit when the NUT charge becomes infinity, the compact dimension opens up, resulting in an ALE space. For the $n$-centre solution, if we insist to keep the positions $z_i$ and charges $\alpha_i$ of all centres fixed, in view of (\ref{NUT_charge_infinity}), there are two natural ALE limits, corresponding to $\beta\rightarrow0$ and $\beta\rightarrow+\infty$ respectively. These limits shall be called D-to-D and D-to-O ALE limits:
\begin{align}
	\text{D-to-D ALE limit: }&\quad\beta\rightarrow 0,\text{ with } n_{\text{nut}}\rightarrow +\infty,\\
	\text{D-to-O ALE limit: }&\quad\beta\rightarrow+\infty,\text{ with } n_{\text{nut}}\rightarrow -\infty,
\end{align}
for reasons that will be explained below. The $n$-centre solution thus interpolates between its D-to-D and D-to-O ALE limits, with $\beta$ as the controlling parameter.

\subsubsection{D-to-D ALE limit}

The D-to-D ALE limit of the $n$-centre solution is obtained by directly setting
\begin{equation}
	\beta= 0,
\end{equation}
in Eqs.~(\ref{w_of_multisoliton})--(\ref{B_function}). In view of (\ref{Psi2_with_k2}), the self-dual Weyl tensor has a well-defined non-zero value in this limit, so the resulting solution is still type D on the self-dual side. Recall that the $n$-centre solution has Petrov type D$^+$X$^-$ for some $X\in\{\text{I, D, O}\}$. Its D-to-D ALE limit will have Petrov type D$^+$Y$^-$, with Y more special than or the same as X. {We remark that the D-to-D ALE limit identified here is equivalent to the Hermitian ALE metrics studied by Araneda and Lucietti \cite{Araneda:2025uqo}.}

Here, we confirm the asymptotic structure of the $n$-centre solution in the D-to-D ALE limit. We first define coordinates $(r,\theta)$ by the equations
\begin{equation}
	\label{CT_ALE_limit}
	\rho=\frac{1}{4}r^2\sin\theta,\quad z=\frac{1}{4}r^2\cos\theta.
\end{equation}
The infinity is located at $r\rightarrow\infty$. Near infinity, at large $r$, the various functions (\ref{w_rho})--(\ref{B_function}) that appear in the metric (\ref{metric_onesided_typeD}) behave as follows:
\begin{equation}
	\begin{aligned}
		w_\rho&\approx-\sin\theta\sum\alpha_i-\frac{4\sin\theta\cos\theta}{r^2}\sum\alpha_iz_i-\frac{8\sin\theta(3\cos^2\theta-1)}{r^4}\sum\alpha_iz_i^2,\\
		w_z &\approx -\cos\theta\sum\alpha_i+\frac{4\sin^2\theta}{r^2}\sum\alpha_iz_i+\frac{24\sin^2\theta\cos\theta}{r^4}\sum \alpha_iz_i^2,\\
		w&\approx-\frac{r^2}{4}\sum\alpha_i+\cos\theta\sum\alpha_i z_i-\frac{2\sin^2\theta}{r^2}\sum \alpha_iz_i^2,\\
		L&\approx\left(\sum\alpha_i\right)^2-\frac{16\sin^2\theta}{r^4}\left[\sum\alpha_i\sum\alpha_iz_i^2-\left(\sum\alpha_i z_i\right)^2\right],\\
		B&\approx\frac{r^4\cos\theta}{16}\sum\alpha_i-\frac{r^2(1+\cos^2\theta)}{4}\sum\alpha_iz_i+\frac{\cos\theta(3-\cos^2\theta)}{2}\sum \alpha_iz_i^2,
	\end{aligned}
\end{equation}
where higher-order terms have been omitted.

Under the coordinate transformation (\ref{CT_ALE_limit}), we have $\dif\rho^2+\dif z^2=\frac{r^2}{4}\left(\dif r^2+\frac{r^2}{4}\dif\theta^2\right)$. If we choose
\begin{equation}
	k_1=\frac{\sum\alpha_i\sum\alpha_iz_i^2-\left(\sum\alpha_i z_i\right)^2}{\sum\alpha_i},\quad k_2=\frac{1}{\sum\alpha_i},
\end{equation}
the $n$-centre solution approaches the following form at large $r$
\begin{equation}
	\dif s^2\approx\dif r^2+\frac{r^2}{4}\left(\dif \tau+\cos\theta\dif \phi\right)^2+\frac{r^2}{4}\left(\dif\theta^2+\sin^2\theta\dif\phi^2\right).
\end{equation}
Higher-order fall-off behaviour of the metric-tensor components near infinity can be determined, but we will not pursue here. The D-to-D ALE limit of the $n$-centre solution is seen to be ALE.

\subsubsection{D-to-O ALE limit}

\label{sec_D_to_O_limit}

The $n$-centre solution admits a degenerate limit in which the self-dual Weyl tensor vanishes. This limit is obtained by first setting $k_1$ and $k_2$ as
\begin{equation}
	k_1=\beta^2,\quad k_2=1,
\end{equation}
respectively, and then taking the limit
\begin{equation}
	\beta\rightarrow \infty.
\end{equation}
We notice that in this limit, 
\begin{equation}
	w\approx -\beta,
\end{equation}
to leading order, so the Weyl scalar $\Psi_2$ as given by (\ref{Psi2_with_k2}) vanishes: 
\begin{equation}
	\label{Psi2_D-to-O-limit}
	\Psi_2\rightarrow 0.
\end{equation}
The self-dual Weyl tensor in the resulting solution thus degenerates to type O (see Eq.~(\ref{self-dual-Weyl-tensor})). If the $n$-centre solution has Petrov type D$^+$X$^-$, its D-to-O ALE limit will have Petrov type O$^+$Y$^-$. It turns out that Y is always equal to X (see Table~\ref{table_multisoliton}).

We now determine the resulting metric of the $n$-centre solution in the D-to-O ALE limit explicitly. The function  $L$ is unchanged in this limit. The function $W$ and the one-form $A$ become
\begin{align}
	W&\rightarrow VL^{-1}, \quad V\equiv -\rho^{-1}w_\rho,\\
	A&\rightarrow w_z L^{-1}\dif \phi=A_\phi\dif\phi.
\end{align}
We can calculate $g_{\phi\phi}$ as
\begin{align}
	g_{\phi\phi}&=W^{-1}A_\phi^2+\rho^2W=V^{-1}L^{-1}(w_z^2+\rho^2V^2)=V^{-1}.
\end{align}
The metric (\ref{metric_onesided_typeD}) in this limit can then be written as
\begin{equation}
	\label{metric_D-to-O-limit}
	\dif s^2=V^{-1}\left(\dif\phi+F\dif \tau\right)^2+V(\rho^2\dif \tau^2+\dif \rho^2+\dif z^2),
\end{equation}
where the functions $V$ and $F$ are given explicitly by
\begin{align}
	\label{functions_in_GH}
	V=-x_z=\sum\frac{\alpha_i}{\sqrt{\rho^2+(z-z_i)^2}},\quad	F=w_z=-\sum\fracb{\alpha_i(z-z_i)}{\sqrt{\rho^2+(z-z_i)^2}}.
\end{align}
We recognise (\ref{metric_D-to-O-limit})--(\ref{functions_in_GH}) as the multi-centre Gibbons--Hawking solutions \cite{Gibbons:1978tef}. Since both the Ricci tensor and the self-dual Weyl tensor vanish, the limiting metric (\ref{metric_D-to-O-limit}) is hyper-K{\"a}hler. We can check directly that the fundamental two-form (\ref{Khaler_form}) in the D-to-O limit becomes
\begin{equation}
	\omega=(\dif\phi+F\dif\tau)\wedge\dif z+\rho V\dif\rho\wedge\dif\tau,
\end{equation}
which is closed $\dif\omega=0$. The associated Hermitian structure is thus K{\"a}hler. The ALE asymptotic structure of the multi-centre Gibbons--Hawking solutions obtained above is well known,\footnote{The multi-centre Gibbons--Hawking metrics, with all centres collinearly aligned, are characterised by a function $V=\epsilon+\sum\frac{2n_i}{\sqrt{\rho^2+(z-z_i)^2}}$. The multi-centre Gibbons--Hawking solutions obtained here are a special case of the multi-centre Gibbons--Hawking metrics with $\epsilon=0$. They are ALE and can be regarded as the multi-centre generalisation of the Eguchi--Hanson solution with arbitrary charges.\label{footnote_GH}} so we will not show explicitly here.

We remark that the 1-centre solution has a D-to-O ALE limit with $n_{\text{nut}}\rightarrow -\infty$, but it does not have a D-to-D ALE limit with $n_{\text{nut}}\rightarrow +\infty$. The 2-centre solution has both D-to-D and D-to-O ALE limits, and they result in the same Riemannian space but with different complex structures that determine opposite orientations. For $n\ge 3$, the $n$-centre solution has both D-to-D and D-to-O ALE limits, with distinct resulting ALE spaces (not isometric to each other). More will be said about these limits in the next section when we discuss explicit examples. In Sec.~\ref{sec_ISM_classification}, the D-to-D and D-to-O ALE limits are understood in a unified way when the multi-centre solutions are fitted into the scheme of multi-soliton solutions on flat space.

\subsection{Further regularity considerations}

In this paper, our primary interest is in the local metrics of the multi-centre solutions, with ranges of parameters given by (\ref{range_alpha}) and (\ref{range_beta}). These solutions are ``locally regular'' in the sense as described in \cite{Chen:2015iex}. They are well-behaved on the half-plane with $(\rho>0,-\infty<z<\infty)$, but may contain mild conical and orbifold singularities along the $z$-axis with $\rho=0$. If the interest is to find completely regular Ricci-flat solutions, i.e., gravitational instantons, further conditions need to be imposed on the local metrics that we have obtained. One such necessary condition is studied in appendix \ref{sec_proof_nge3}.

\section{Examples}

\label{sec_examples}

In this section, we study the multi-centre solutions with number of centres $n\le 3$ in detail. Many closed-form one-sided type-D Ricci-flat metrics are recovered from them.

\subsection{The 1-centre solution}

\subsubsection{1-centre Taub--NUT}

The 1-centre solution is determined by the following pair of generating potentials
\begin{equation}
	w=-\beta-\alpha_1\sqrt{\rho^2+(z-z_1)^2}=-\beta-\frac{\alpha_1(\rho^2+\mu_1^2)}{2\mu_1},\quad x=\alpha_1 \ln \frac{\mu_1}{\rho} ,
\end{equation}
and it turns out to be the 1-centre Taub--NUT solution \cite{Newman:1963yy,Hawking:1976jb}. With the constants $k_1$, $k_2$ and $z_1$ chosen as 
\begin{equation}
	k_1=\beta, \quad k_2=\alpha_1^{-1}, \quad z_1=0,
\end{equation}
 respectively, the 1-centre solution can be written down explicitly as
\begin{align}
	\dif s^2&=W^{-1}\left(\dif \tau-\frac{\beta }{\alpha_1}\frac{z}{\sqrt{\rho^2+z^2}}\dif\phi\right)^2+W(\rho^2\dif\phi^2+\dif\rho^2+\dif z^2),\\ W&=1+\frac{\beta}{\alpha_1}\frac{1}{\sqrt{\rho^2+z^2}}.\nonumber
\end{align}
If we perform coordinate transformations and redefine parameter as follows
\begin{equation}
	\rho=r\sin\theta,\quad z=r\cos\theta,\quad 2n=\frac{\beta}{\alpha_1},
\end{equation}
the 1-centre solution is brought to the 1-centre Taub--NUT solution in the familiar form \cite{Hawking:1976jb} 
\begin{align}
	\label{self-dual Taub-NUT}
	\dif s^{2}=W^{-1} { \left( {\dif\tau}-2n\cos  \theta\,
		{\dif\phi} \right) ^{2}}+W \left( {\dif r}^{2}+{r}^{2}{
		{\dif\theta}}^{2}+{r}^{2} \sin ^2 \theta \,{{\dif\phi}}^{2} \right),\quad W=1+\frac{2 n}{r}\,.
\end{align}
The asymptotic NUT charge $n_{\text{nut}}$ (cf. (\ref{NUT-charge}) and (\ref{metric_ALF})) is identified as
\begin{equation}
	n_{\text{nut}}=-n<0,
\end{equation}
which is always negative. We remark that the 1-centre Taub--NUT solution was first recovered from the 1-centre solution in equivalent forms in \cite{Biquard:2021gwj}.

The fundamental two-form (\ref{Khaler_form}) associated with the Hermitian structure arising from the one-sided type-D condition is calculated as
\begin{equation}
	\omega=-(\dif\tau-2n\cos\theta\dif\phi)\wedge\dif r+r(r+2n)\sin\theta\dif\theta\wedge\dif\phi.
\end{equation}
The 1-centre solution, which has been shown to be the 1-centre Taub--NUT solution, has Petrov type D$^+$O$^-$. It does not admit a D-to-D ALE limit: taking the limit $\beta=0$ results in flat space metric, which is no longer of Petrov type D on the self-dual side. On the other hand, the D-to-O ALE limit exists: taking the limit $\beta\rightarrow \infty$ results in flat space with Petrov type O$^+$O$^-$.

\subsection{The 2-centre solution}

The 2-centre solution is determined by the following pair of generating potentials
\begin{align}
	\label{generating_function_2centre}
	w=-\beta-\sum\limits_{i=1}^{2}\alpha_i\sqrt{\rho^2+(z-z_i)^2}=-\beta-\sum\limits_{i=1}^{2}\frac{\alpha_i(\rho^2+\mu_i^2)}{2\mu_i},\quad x=\sum\limits_{i=1}^{2}\alpha_i \ln \frac{\mu_i}{\rho},
\end{align}
and it turns out to be the Kerr--NUT solution \cite{Carter:1968ks,Gibbons:1979nf}, as expected from the combined properties that it possesses: (a) two turning points, (b) Petrov type D on one side and (c) ALF asymptotic structure.

\subsubsection{Kerr--NUT}
To transform the 2-centre solution to the familiar form of the Kerr--NUT solution, we first redefine parameters
\begin{equation}
	z_1=-l,\quad z_2=l,
\end{equation}
and introduce prolate spheroidal coordinates $(p,q)$ by their relation to the Weyl--Papapetrou coordinates $(\rho,z)$
\begin{equation}
	\rho=l\sqrt{(1-p^2)(q^2-1)},\quad z=lpq.
\end{equation}
A major advantage of using prolate spheroidal coordinates is that now the solitons $\mu_1$ and $\mu_2$ become rational functions (polynomials in fact) of $(p,q)$ 
\begin{equation}
	\mu_1=l(1-p)(q-1),\quad \mu_2=l(1-p)(q+1),
\end{equation}
and so do all metric-tensor components. The full conversion of the 2-centre solution to prolate spheroidal coordinates can be handled much faster using a computer algebra system. We do not record all metric-tensor components in $(p,q)$ coordinates here.

Next, we define Boyer--Lindquist coordinates $(r,\theta)$ by the equations
\begin{align}
	\label{rtheta1topq}
	r=l q+\frac{4\alpha_1\alpha_2 l^2+\beta^2}{2\beta(\alpha_1+\alpha_2)}\,,\quad \cos\theta=p\,,
\end{align}
and redefine coordinate $\tau$ by doing the substitution 
\begin{equation}
	\tau\rightarrow\tau+\frac{2l(\alpha_1-\alpha_2)}{\alpha_1+\alpha_2}\phi\,.
\end{equation}
After that, we fix the parameters $k_1$ and $k_2$ as
\begin{equation}
	k_1=\beta,\quad k_2=(\alpha_1+\alpha_2)^{-1}.
\end{equation}
Then it is straightforward to show that the 2-centre solution is brought to the Kerr--NUT solution in the following form
\begin{align}
	\label{metric_Kerr-NUT}
	{\dif s}^{2}=&{\frac {{\Delta_2}}{{\Sigma}}} ({\dif\tau}+ A) ^{2}+{\frac {\Sigma\Delta_1
		}{{\Delta_2}}}\sin ^{2} \theta\,\dif\phi ^{2}+{\Sigma} \left( {\frac {{\dif r}^{2}}{{\Delta_1}}}
	+{{\dif\theta}}^{2} \right),\\
	{\Sigma}=&{r}^{2}- \left( n-a\cos \theta  \right) ^{
		2},\quad
	{\Delta_1}={r}^{2}-2mr-{a}^{2}+{n}^{2}\,,\quad {\Delta_2}=\Delta_1+a^2\sin^2\theta,\nonumber\\
	A=&\left( 2n
	\cos \theta-\frac{2(mr+an\cos\theta-n^2)}{\Delta_2}a \sin  ^{2} \theta
	\right) {\dif\phi},\nonumber
\end{align}
with the parameters $m$, $n$, and $a$ identified as, respectively,
\begin{align}\label{parameters_2center_KN}  m=\frac{4\alpha_1\alpha_2 l^2+\beta^2}{2\beta(\alpha_1+\alpha_2)}\,,\quad n=\frac{4\alpha_1\alpha_2 l^2-\beta^2}{2\beta(\alpha_1+\alpha_2)}\,,\quad a=\frac{l(\alpha_1-\alpha_2 )}{\alpha_1+\alpha_2}\,.
\end{align}

The fundamental two-form (\ref{Khaler_form}) associated with the Hermitian structure arising from the one-sided type-D condition can be calculated in Boyer--Lindquist coordinates as
\begin{equation}
	\omega=-\dif\tau\wedge(\dif r-a\sin\theta\dif\theta)+\left[(r^2-a^2-n^2)\sin\theta\dif\theta+(2n\cos\theta+a\sin^2\theta)\dif r\right]\wedge\dif\phi.
\end{equation}
It is well known that the Kerr--NUT solution has Petrov type D$^+$D$^-$, so does the 2-centre solution.

{The Taub-bolt instanton \cite{Page:1978hdy} arises as a special case of the Kerr-NUT solution with $(a=0,n=\pm \frac{4}{5}m)$ (see e.g., \cite{Chen:2010zu}). In view of Eq.~(\ref{parameters_2center_KN}) and the identification of the Kerr-NUT solution with the 2-centre solution, the Taub-bolt instanton is seen to arise as a special case of the 2-centre solution with parameters given by
\begin{equation}
	\alpha_2=\alpha_1=\frac{3\beta}{2l}\text{ or }\frac{\beta}{6l}.
\end{equation} 
The mass parameter of the Taub-bolt instanton is then identified as
\begin{equation}
	m=\frac{5l}{3}.
	\end{equation}}

\subsubsection{2-centre Gibbons--Hawking}

The D-to-D ALE limit of the 2-centre solution is obtained by setting $\beta=0$ in the generating potential $w$ given in Eq.~(\ref{generating_function_2centre}). We show that it results in the 2-centre Gibbons--Hawking solution \cite{Gibbons:1978tef}. 

First, we reparametrise the centre positions $z_1$ and $z_2$ as
\begin{equation}
	z_1=-l,\quad z_2=l, 
\end{equation}
respectively. Then, we introduce coordinates $(r,\theta)$ by
\begin{equation}
	\rho=r\sin\theta, \quad z=r\cos\theta.
\end{equation}
With the constants $k_1$ and $k_2$ chosen as
\begin{equation}
	k_1=8l^2,\quad k_2=\frac{2}{\alpha_1\alpha_2},
\end{equation}
respectively, the 2-centre solution becomes
\begin{equation}
	\label{double-centered Gibbons-Hawking}
	{\dif s}^{2}=V^{-1} { \left( {\dif\phi}+A \right) ^{2}}+V
	( {\dif r}^{2}+{r}^{2} {{\dif\theta}}^{2}+{r}^{2} \sin^{2}
	\theta\, {{\dif\tau}}^{2} )\,,
\end{equation}
where the function $V$ and the one-form $A$ are given by
\begin{align}
	\label{metric_2GH}
	V&={\frac {2 n_1}{{r_1}}}+{\frac {2 n_2}{{r_2}}}\,,\quad r_1=\sqrt {{r}^{2}+{l}^{2}+2lr\cos  \theta  },\quad r_2=\sqrt {{r}^{2}+{l}^{2}-2lr\cos  \theta  },\cr
	A&={\frac {2 n_1 \left( r\cos \theta  +l \right) }{{r_1}}}\,{\dif\tau
	}+{\frac {2 n_2 \left( r\cos \theta -l
			\right) }{{r_2}}}\,{\dif\tau}\,,
\end{align}
with parameters $n_1$ and $n_2$ identified as
\begin{equation}
	n_1=\frac{1}{\alpha_1},\quad n_2=\frac{1}{\alpha_2},
\end{equation}
respectively. We see that this is the 2-centre Gibbons--Hawking solution. 

The fundamental two-form (\ref{Khaler_form}) associated with the Hermitian structure arising from the one-sided type-D condition is calculated as
\begin{equation}
	\label{fundamental_twoform_2centre}
	\omega=\frac{1}{2l}\dif\phi\wedge(r_2\dif r_1-r_1\dif r_2)+2(n_2\dif r_1+n_1\dif r_2)\wedge \dif\tau,
\end{equation}
which can be verified directly to be non-closed. The 2-centre Gibbons--Hawking solution, with the fundamental two-form given by (\ref{fundamental_twoform_2centre}), has Petrov type D$^+$O$^-$. So, in the current case for the 2-centre solution, the D-to-D ALE limit described above converts Petrov type from D$^+$D$^-$ to D$^+$O$^-$. The D-to-O ALE limit is obtained by taking $\beta\rightarrow \infty$ as in the general case described in Sec.~\ref{sec_D_to_O_limit}, with a resulting solution that is also the 2-centre Gibbons--Hawking solution. This limit will convert the Petrov type of the self-dual Weyl tensor from D to O. The anti-self-dual Weyl tensor retains its Petrov type D. So, for the 2-centre solution, the D-to-O ALE limit converts its Petrov type from D$^+$D$^-$ to O$^+$D$^-$.

\subsubsection{Eguchi--Hanson}

The Eguchi-Hanson solution \cite{Eguchi:1978xp} is a special case of the 2-centre Gibbons--Hawking solution with the two centres carrying equal charges ($n_1=n_2$ in (\ref{metric_2GH})). Similar to the 2-centre Gibbons--Hawking solution, it can be obtained both from the D-to-D and D-to-O ALE limits from the 2-centre solution. Here, we will only describe how the Eguchi--Hanson solution with Petrov type D$^+$O$^-$ arises from the D-to-D limit of the 2-centre solution.

First, we set parameters $\beta$ and $\alpha_2$ in the generating potential $w$ of the 2-centre solution given in (\ref{generating_function_2centre}) as
\begin{equation}
	\beta=0,\quad \alpha_2=\alpha_1, 
\end{equation}
and reparametrise the centre positions $z_1$ and $z_2$ as
\begin{equation}
	z_1=-\frac{a^2}{4},\quad z_2=\frac{a^2}{4},
\end{equation}
respectively. Then, we define coordinates $(r,\theta)$ by
\begin{equation}
	\rho=\frac{1}{4}\sqrt{r^4-a^4}\sin\theta,\quad z=\frac{r^2}{4}\cos\theta.
\end{equation}
With the constants $k_1$ and $k_2$ chosen as
\begin{equation}
	k_1=\frac{\alpha_1 a^4}{8},\quad k_2=\frac{1}{2\alpha_1},
\end{equation}
respectively, the 2-centre solution is brought to the form 
\begin{align}
	\label{Eguchi-Hanson}
	{\dif s}^{2}=\left( 1-{\frac {{a}^{4}}{{r}^{4}}} \right) \frac{{r}^{
			2}}{4} \left( {\dif\tau}+\cos \theta\, {\dif\phi} \right) ^{2
	}+\left( 1-{\frac {{a}^{4}}{{r}^{4}}} \right) ^{-1} {\dif r}^{2}+\frac{{r}^{2}}{4} ( {{\dif\theta}}^{2}+ \sin^{2} \theta\,
	{{\dif\phi}}^{2} )\,,
\end{align}
which is recognised as the Eguchi--Hanson solution. We remark that the Eguchi--Hanson solution was first recovered from the 2-centre solution in equivalent forms in \cite{Biquard:2023rbr,Araneda:2025uqo}. The fundamental two-form (\ref{Khaler_form}) associated with the Hermitian structure arising from the one-sided type-D condition is calculated as
\begin{equation}
	\omega=-\frac{r}{2}(\dif\tau+\cos\theta\dif\phi)\wedge\dif r+\frac{r^2}{4}\sin\theta\dif\theta\wedge\dif\phi.
\end{equation}

\subsection{The 3-centre solution}

The 3-centre solution is determined by the following pair of generating potentials
\begin{align}
	\label{generating_function_3centre}
	w=-\beta-\sum\limits_{i=1}^{3}\alpha_i\sqrt{\rho^2+(z-z_i)^2}=-\beta-\sum\limits_{i=1}^{3}\frac{\alpha_i(\rho^2+\mu_i^2)}{2\mu_i},\quad x=\sum\limits_{i=1}^{3}\alpha_i \ln \frac{\mu_i}{\rho} ,
\end{align}
and it turns out to be the Chen--Teo solution \cite{Chen:2015vva,Biquard:2023rbr}. 

\subsubsection{Chen--Teo}

The Chen--Teo solution was written in C-metric-like coordinates $(x,y)$,\footnote{The C-metric-like coordinates $(x,y)$ introduced here are distinct from the $(x,y)$ coordinates used in the LeBrun--Tod ansatz (\ref{metric_Lebrun_Tod}). The C-metric-like coordinate $x$ is also distinct from the generating potential $x$. We are currently running out of symbols for coordinates.} which are related to the Weyl--Papapetrou coordinates by the following equations
\begin{equation}
	\label{C-metric-like-coordinates}
	\rho=\frac{\sqrt{-XY}}{(x-y)^2}, \quad \dif\rho^2+\dif z^2\propto \frac{\dif x^2}{X}-\frac{\dif y^2}{Y},
\end{equation}
where $X=X(x)$ and $Y=Y(y)$, called structure functions of the solution, are polynomials of their respective arguments and satisfy $Y(x)=X(x)$. The metrics on the two sides of the symbol $\propto$ are conformal to each other. Consider the inverse coordinate transformation from $(x,y)$ to $(\rho,z)$. The integrability condition for the coordinate $z$ derived from Eq.~(\ref{C-metric-like-coordinates}) is satisfied if $X$ and $Y$ are at most quartic (see Appendix H of \cite{Harmark:2004rm}). The static C-metric (see, e.g., \cite{Emparan:2001wk}) has cubic structure functions. The Chen--Teo solution has quartic structure functions \cite{Chen:2015vva}. In \cite{Chen:2015iex}, the Chen--Teo solution was re-derived by first using C-metric-like coordinates with cubic structure functions, and then applying a M\"{o}bius transformation of the form $(x\rightarrow\frac{ax+b}{cx+d},y\rightarrow\frac{ay+b}{cy+d})$, which turns the structure functions to quartic ones. The form of Eq.~(\ref{C-metric-like-coordinates}) is preserved under the M\"{o}bius transformation.

The basic idea of using C-metric-like coordinates to simplify the 3-centre solution is to absorb as many parameters as possible into the structure functions. Here, we will use C-metric-like coordinates with quartic structure functions directly, and explicitly show that the 3-centre solution is the Chen--Teo solution. We begin with the following quartic structure function
\begin{equation}
	\label{X}
	X=a_0+a_1x+a_2x^2+a_3x^3+a_4x^4=a_4(x-x_1)(x-x_2)(x-x_3)(x-x_4).
\end{equation}
The two sets of parameters $(a_0,a_1,a_2,a_3,a_4)$ and $(x_1,x_2,x_3,x_4,a_4)$ are related to each other and the relation is determined from the above equation. Now treat the roots $(x_1,x_2,x_3)$ as fixed parameters whose values can be arbitrarily assigned. We can determine $a_4$ and the remaining root $x_{4}$ from the centre positions $(z_1,z_2,z_3)$ by imposing the following two equations
\begin{equation}
	\label{3soliton_rod_length_parameters}
	z_2-z_1=\frac{a_4(x_1-x_4)(x_2-x_3)}{2}, \quad z_3-z_2=\frac{a_4(x_1-x_2)(x_3-x_4)}{2}.
\end{equation}
We then exploit the translational symmetry (\ref{symmetry_translational}) of the solution, and express the centre positions $(z_1,z_2,z_3)$ as
\begin{equation}
	\label{3soliton_positions}
	z_1=-\frac{a_4(x_1x_2+x_3x_4)}{2},\quad z_2=-\frac{a_4(x_1x_3+x_2x_4)}{2},\quad z_3=-\frac{a_4(x_1x_4+x_2x_3)}{2}.
\end{equation}
We change perspective and take $(x_4,a_4)$ as independent parameters that encode the centre positions $(z_1,z_2,z_3)$. With the structure function given in (\ref{X}), the coordinate $z$ can be integrated in terms of $(x,y)$ as
\begin{equation}
\label{3soliton_WP}
 z=\frac{2(a_0+a_2xy+a_4x^2y^2)+(x+y)(a_1+a_3xy)}{2(x-y)^2}.
\end{equation}
The integration constant has been fixed so that $z(x=x_2,y=x_1)=z_1$, where $z_1$ is given in (\ref{3soliton_positions}).

Recall that a soliton $\mu_i=\sqrt{\rho^2+(z-z_i)}-(z-z_i)$ contains a square-root expression in terms of coordinates $(\rho,z)$ or $(x,y)$. The centre positions $(z_1,z_2,z_3)$, given by (\ref{3soliton_positions}), are defined in such a way that square roots are eliminated for all the three solitons ($\mu_1,\mu_2,\mu_3$) in $(x,y)$ coordinates: the expressions in all square roots become perfect squares. More specifically, now $(\mu_1,\mu_2,\mu_3)$ attain simple rational forms and they can be written as
\begin{equation}\label{3soliton_solitons}
	\begin{aligned}
		\mu_1&=-\frac{a_4(x-x_1)(x-x_2)(y-x_3)(y-x_4)}{(x-y)^2},\\
		\mu_2&=-\frac{a_4(x-x_1)(x-x_3)(y-x_2)(y-x_4)}{(x-y)^2},\\
		\mu_3&=-\frac{a_4(x-x_1)(x-x_4)(y-x_2)(y-x_3)}{(x-y)^2}.
	\end{aligned}
\end{equation}
The expression following the second equal sign in the generating potential 
$w$ of Eq.~(\ref{generating_function_3centre}) is used in our calculations. The expression following the first equal sign can also be used, and all square roots in it can be eliminated as well. All metric-tensor components can be converted in terms of coordinates $(x,y)$, and they are all rational functions. This is a major step in simplifying the 3-centre solution.

At this point, the 3-centre solution has parameters $(\beta,\alpha_1,\alpha_2,\alpha_3,x_1,x_2,x_3,x_{4},a_4)$ in addition to $(k_1,k_2)$. Now we redefine $\beta$ in terms of a new parameter $\nu$ as
\begin{equation}
	\label{beta_CT}
	\beta=\frac{1-\nu}{2(1+\nu)}.
\end{equation}
Recall that we have the freedom to assign arbitrary values to the root parameters $(x_1,x_2,x_3)$. By using this freedom, we can simplify the generating potential $w$ so that it takes a simple form $w\propto\frac{\nu x+y}{x-y}$, by imposing the following equations
\begin{equation}
	\label{alpha_CT}
	\begin{aligned}
		\alpha_1&=\frac{x_3x_4-x_1x_2}{a_4(x_1-x_3)(x_1-x_4)(x_2-x_3)(x_2-x_4)},\\
		\alpha_2&=\frac{x_2x_4-x_1x_3}{a_4(x_1-x_2)(x_1-x_4)(x_3-x_2)(x_3-x_4)},\\
		\alpha_3&=\frac{x_2x_3-x_1x_4}{a_4(x_1-x_2)(x_1-x_3)(x_4-x_2)(x_4-x_3)}.
	\end{aligned}
\end{equation}
Now we take $(x_1,x_2,x_3)$ instead of $(\alpha_1,\alpha_2,\alpha_3)$ as the independent parameters. To further simplify, we choose parameters $k_1$ and $k_2$ as 
\begin{equation}
	k_1=\frac{1}{2(1+\nu)^2},\quad k_2=\sqrt{k},
\end{equation}
respectively, and redefine coordinates by doing the substitution
\begin{equation}
	\tau\rightarrow \tau+\frac{1}{2}a_2(1+2\nu-\nu^2)\phi,
\end{equation}
followed by
\begin{equation}
	\phi\rightarrow -\phi.
\end{equation}
After changing parameters from $(x_1,x_2,x_3,x_4,a_4)$ to $(a_0,a_1,a_2,a_3,a_4)$, the 3-centre solution is brought to the Chen--Teo solution given in \cite{Chen:2015vva}
\begin{align}
	\label{metric_new_ricci}
	\dif s^2&=\frac{\left(F\dif \tau+G\dif \phi\right)^2}{(x-y)HF}+\frac{kH}{(x-y)^3}\left(\frac{\dif x^2}{X}-\frac{\dif y^2}{Y}-\frac{XY}{kF}\,\dif \phi^2\right),\\
	H&=(\nu x+ y)[(\nu x-y)(a_1-a_3xy)-2(1-\nu)(a_0-a_4x^2y^2)]\,,\quad F=y^2X-x^2Y\,,\nonumber\\
	G&=(\nu^2a_0+2\nu a_3y^3+2\nu a_4y^4-a_4y^4)X+(a_0-2\nu a_0-2\nu a_1x-\nu^2 a_4x^4)Y\,,\nonumber\\
	X&=a_0+a_1x+a_2x^2+a_3x^3+a_4x^4\,,\quad Y=a_0+a_1y+a_2y^2+a_3y^3+a_4y^4\,.\nonumber
\end{align}
We remark that the generating potential $w$ attains a simple form after all the above transformations:
\begin{equation}
	w=\frac{\nu x+y}{(1+\nu)(x-y)}.
\end{equation}
This expression, according to our general analysis in Sec.~\ref{sec_one_sided_type_D_metrics}, completely determines the 3-centre solution. The information of the centre positions $(z_1,z_2,z_3)$ and the centre charges $(\alpha_1,\alpha_2,\alpha_3)$ are all encoded in the C-metric-like coordinate functions ($x,y$).

The fundamental two-form (\ref{Khaler_form}) associated with the Hermitian structure arising from the one-sided type-D condition is calculated as
\begin{equation}
	\omega=\frac{\sqrt{k}}{(x-y)^2F}\left[\left(F\dif\tau+G\dif\phi\right)\wedge(-y\dif x+x\dif y)+\frac{HXY}{x-y}\left(\frac{x\dif x}{X}-\frac{y\dif y}{Y}\right)\wedge\dif\phi\right],
\end{equation}
in agreement with results found in the work of Aksteiner and Anderson \cite{Aksteiner:2021fae}.
The 3-centre solution has Petrov type D$^+$I$^-$ by the same work. It admits both D-to-D and D-to-O ALE limits. The D-to-D ALE limit corresponds to $\beta\rightarrow 0$ or equivalently $\nu\rightarrow 1$ in view of the expression of $\beta$ in (\ref{beta_CT}). The resulting solution is the Pleba\'nski--Demia\'nski solution with Petrov type D$^+$D$^-$. The D-to-O ALE limit corresponds to $\beta\rightarrow +\infty$ or equivalently $\nu\rightarrow -1$ in view of (\ref{beta_CT}). The resulting solution is the 3-centre Gibbons--Hawking solution with Petrov type O$^+$I$^-$. Both limits were observed in \cite{Chen:2015vva}. Here, we will only discuss the D-to-D ALE limit with a resulting Pleba\'nski--Demia\'nski solution. 

\subsubsection{Pleba\'nski--Demia\'nski}
The D-to-D ALE limit of the 3-centre solution is obtained by setting $\beta=0$ in the generating potential $w$ given in Eq.~(\ref{generating_function_3centre}). We will follow similar steps as described in the previous subsection to bring the resulting solution to the familiar form of the Pleba\'nski--Demia\'nski solution \cite{Plebanski:1976gy}. 

We use C-metric-like coordinates defined by Eq.~(\ref{C-metric-like-coordinates}) with the following structure function
\begin{equation}
	\label{PD_structure_function}
	X=\gamma+2nx-\epsilon x^2+2mx^3+\gamma x^4=\gamma(x-x_1)(x-x_2)(x-x_3)(x-x_4),
\end{equation}
with a constraint on the four roots
\begin{equation}
	\label{PD_roots}
	x_1x_2x_3x_4=1.
\end{equation}
The two sets of parameters $(m,n,\epsilon,\gamma)$ and $(x_1,x_2,x_3,\gamma)$ are related to each other and their relation is determined from the equation (\ref{PD_structure_function}). This is the same  situation as in the previous subsection with $a_4=\gamma$ and the extra condition (\ref{PD_roots}). Taking $(x_1,x_2)$ as fixed parameters whose values can be arbitrarily assigned, by imposing (\ref{3soliton_rod_length_parameters}), we can determine $\gamma$ and $x_3$ from rod lengths $z_2-z_1$ and $z_3-z_2$ of the 3-centre solution. By exploiting the translational symmetry (\ref{symmetry_translational}), we reparametrise the centre positions $(z_1,z_2,z_3)$ by the equation (\ref{3soliton_positions}). The coordinate $z$ and three solitons $(\mu_1,\mu_2,\mu_3)$ in terms of $(x,y)$ are given by equations (\ref{3soliton_WP}) and (\ref{3soliton_solitons}) respectively. The metric-tensor components are then converted in terms of $(x,y)$ coordinates and they all become rational functions.

At this point, the solution has parameters $( x_1,x_2,x_3,\gamma,\alpha_1,\alpha_2,\alpha_3)$ in addition to $k_1$ and $k_2$. We have the freedom to assign arbitrary values to the parameters $(x_1, x_2)$. Now by using this freedom and exploiting the parameter-rescaling symmetry (\ref{symmetry_scaling}), we can impose the following conditions without loss of generality
\begin{equation}
	\begin{aligned}
		\label{alpha_in_PD}
		\alpha_1&=-\frac{2(1+x_1x_2)}{(x_1-x_3)(x_1-x_4)(x_2-x_3)(x_2-x_4)(1-x_1x_2)},\\
		\alpha_2&=-\frac{2(1+x_1x_3)}{(x_1-x_2)(x_1-x_4)(x_3-x_2)(x_3-x_4)(1-x_1x_3)},\\
		\alpha_3&=\frac{2(1+x_2x_3)}{(x_1-x_2)(x_1-x_3)(x_4-x_2)(x_4-x_3)(1-x_2x_3)},
	\end{aligned}
\end{equation}
to simplify the expression of the generating potential $w$ so that it has a simple form $w\propto\frac{1-xy}{x-y}$. There is a one-parameter redundancy among the three parameters $(\alpha_1,\alpha_2,\alpha_3)$ due to the parameter-rescaling symmetry, which has been used to make the equations in (\ref{alpha_in_PD}) look balanced between the $\alpha$'s. So we trade $(\alpha_1,\alpha_2,\alpha_3)$ with $(x_1,x_2)$ and take the latter as independent parameters.
Next, we fix $k_1$ and $k_2$ as
\begin{equation}
	k_1=\frac{\gamma^2}{2},\quad k_2=\frac{(m-n)}{\sqrt{2}\gamma^2},
\end{equation}
respectively, and perform coordinate transformations by doing the following simultaneous substitutions
\begin{equation}
	\tau\rightarrow \frac{1}{\sqrt{2}}(\tau+\phi),\quad \phi\rightarrow \frac{1}{\sqrt{2}}(-\tau+\phi).
\end{equation}
After changing parameters from $(x_1,x_2,x_3,\gamma)$ to $(m, n,\epsilon,\gamma)$, the D-to-D ALE limit of the 3-centre solution is brought to the form
\begin{align}
	\label{metric_PD}
	\dif s^2&=\frac{1}{(x-y)^2}\Big[\frac{1-x^2y^2}{X}\dif x^2+\frac{X}{1-x^2y^2}(\dif \phi-y^2\dif \tau)^2\nonumber\\
	&\hspace{0.56in}-\frac{1-x^2y^2}{Y}\dif y^2-\frac{Y}{1-x^2y^2}(\dif \tau-x^2\dif \phi)^2\Big],\\
	X&=\gamma+2nx-\epsilon x^2+2mx^3+\gamma x^4\,,\quad Y=\gamma+2ny-\epsilon y^2+2my^3+\gamma y^4\,,\nonumber
\end{align}
which is recognised as the Pleba\'nski--Demia\'nski solution \cite{Plebanski:1976gy}. The generating potential $w$ attains a simple form after all the above transformations:
\begin{equation}
	w=\frac{\gamma^2 (1-xy)}{(m-n)(x-y)}.
\end{equation}
The fundamental two-form (\ref{Khaler_form}) associated with the Hermitian structure arising from the one-sided type-D condition is calculated as
\begin{equation}
	\omega=\frac{1}{(x-y)^2}\left[(\dif\tau-x^2\dif\phi)\wedge\dif y+\dif x\wedge(\dif\phi-y^2\dif\tau)\right].
\end{equation}

\section{Solitonic classification of the multi-centre solutions}

\label{sec_ISM_classification}

In this section, we discuss how the multi-centre solutions discussed in the previous sections fit into the scheme of multi-soliton solutions on flat space \cite{Chen:2015iex} constructed by the present author using the inverse-scattering method.

\subsection{Multi-soliton solutions on flat space}

The ISM of Belinski and Zakharov (BZ) has been very successfully applied to generate Lorentzian solutions \cite{Belinski:2001ph}, such as the Kerr and multi-Kerr black-hole solutions. More recently, it has been shown that ISM can be used to generate Euclidean solutions as well \cite{Chen:2015iex} by the present author. In \cite{Chen:2015iex}, a class of solutions, called gravitational multi-soliton solutions on flat space, were constructed using the ISM. Flat-space seed solutions were employed, as they were found, to the best of the author's knowledge, to be necessary to ensure that the resulting solution is locally regular and does not contain curvature singularities. An $n$-soliton transformation was applied on the seed solutions. The $n$-soliton solution on flat space obtained is either ALF or ALE. A certain number, say $m$, of solitons can be eliminated in a non-trivial way by appropriately fixing the so-called BZ parameters associated with these solitons. Depending on different ways of fixing the BZ parameters, the resulting solution splits into separate classes, collectively denoted as $[n-m]$-soliton solutions on flat space. The $n$-soliton solution on flat space can then be taken as a special case of $[n-m]$-soliton solutions on flat space with $m=0$.

These $[n-m]$-soliton solutions on flat space can be classified by a pair of integers $(s,s_z)$, where $s$ is called the free-soliton number and $s_z$ the phantom-soliton number. The $m$ solitons that are eliminated no longer correspond to turning points of the solution, so they are called phantom solitons. The remaining $n-m$ solitons in the solution are called free solitons. The free-soliton number $s$ is defined as $s=n-m$. The $m$ phantom solitons further divide into, say $m_1$, positive phantom solitons and, say $m_2$, negative phantom solitons, with $m_1+m_2=m$. The phantom-soliton number $s_z$ is defined as the number of positive phantom solitons minus the number of negative phantom solitons: $s_z=m_1-m_2$. A multi-soliton solution on flat space can be represented by the symbol $[n-m_1-m_2]$, for some integers $n$, $m_1$ and $m_2$. The pair of integers $(s=n-m,s_z=m_1-m_2)$ classifies all multi-soliton solutions on flat space.  All explicitly known Ricci-flat locally regular metrics with two Killing vectors seem to be encompassed in this classification scheme. For example, the Kerr--NUT solution has a free-soliton number $s=2$ and a phantom-soliton number $s_z=0$, the 2-centre Taub--NUT solution has a free-soliton number $s=2$ and a phantom-soliton number $s_z=2$, and the Chen--Teo solution has a free-soliton number $s=3$ and a phantom-soliton number $s_z=1$. The reader is referred to Table~3 of \cite{Chen:2015iex} for more details. Only non-negative values of phantom-soliton number $s_z$ were included in that table in view of some symmetry properties of multi-soliton solutions on flat space.

Here, we reorganise the solitonic classification of multi-soliton solutions on flat space as given in Table~3 of \cite{Chen:2015iex},  according to the pair of integers $(s,s-s_z)$, and present it in Table~\ref{table_multisoliton} below.
In this table, we have allowed negative phantom-soliton number, so $s_z$ can take integer values in the range $-s\le s_z\le s$. A multi-soliton solution on flat space represented by $[n-m_1-m_2]$ with $m_1>m_2$ has a positive phantom-soliton number. We can simply flip the sign of all BZ parameters associated with phantom solitons to change the solution to one with a negative phantom-soliton number, which can be represented by $[n-m_2-m_1]$. With the free-soliton number $s$ fixed, the effect of this operation can be equivalently interpreted as follows: the metric remains unchanged, but a given orientation gets reversed. Multi-soliton solutions on flat space with negative phantom-soliton number $-s\le s_z<0$ are placed in cells with dashed boundary lines in Table~\ref{table_multisoliton}. With these solutions included, in the $i$-th row in Table~\ref{table_multisoliton}, each solution has a mirror image about the cell $(s,s-s_z)=(i,i)$ with the same metric and opposite orientation: The $(s,s-s_z)=(2,3)$ and $(2,1)$ solutions are both 2GH but with opposite orientations; the $(s,s-s_z)=(2,2)$ solution is its own image, so the same KN metric with opposite orientations are included in the same cell; etc.. 

\begin{figure}[!t]
	\begin{center}
		\includegraphics[width=\textwidth]{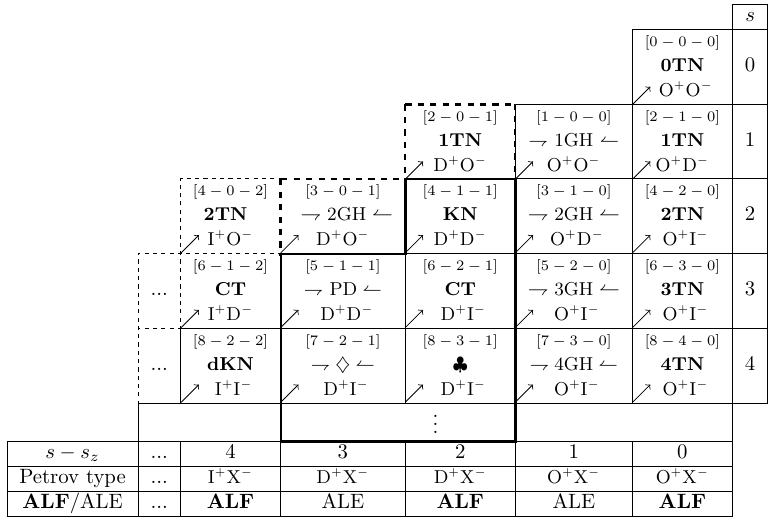}
		\captionsetup{name=Table}
		\caption{Solitonic classification of multi-soliton solutions on flat space constructed in \cite{Chen:2015iex}, arranged according to the pair $(s,s-s_z)$, where $s$ is the free-soliton number and $s_z$ is the phantom-soliton number. For each solution, a solitonic representation of it is given in the form of $[n-m_1-m_2]$ and its Petrov type is indicated below the solution kernel symbol. GH stands for Gibbons--Hawking and TN stands for Taub--NUT; their number of centres is placed in front. So $2$TN stands for the 2-centre Taub--NUT solution. FS stands for flat space, KN stands for Kerr--NUT, PD stands for Pleba\'nski--Demia\'nski, CT stands for Chen--Teo and dKN stands for double Kerr--NUT. The symbols $\diamondsuit$ and $\clubsuit$ denote solutions whose solitonic representations are as indicated. Solutions in boldface, including $\clubsuit$, are ALF; the rest are ALE. The arrows $\leftharpoonup$ and $\rightharpoondown$ represent ALE limits, while the arrows $\nearrow$ represents transition limits.}
		\label{table_multisoliton}
	\end{center}
\end{figure}

We remark that in Table~\ref{table_multisoliton}, the multi-centre Gibbons--Hawking solutions (1GH, 2GH, etc.) refer to the special case of the multi-centre Gibbons--Hawking metrics with $\epsilon=0$ (see footnote~\ref{footnote_GH}), and they are ALE. The multi-centre Taub-NUT solutions (1TN, 2TN, etc.) are another special case with $\epsilon=1$, and they are ALF. It is observed that, in the above table, solutions with even $s-s_z$ are ALF and those with odd $s-s_z$ are ALE.  The arrows $\leftharpoonup$ and $\rightharpoondown$ in the table represent ALE limits, or decompactification limits as called in \cite{Chen:2015iex}. The arrows $\nearrow$ represent transition limits. All three types of limits change the asymptotic structure of a solution.

Petrov types for all explicitly presented solutions ($s\le 4$) in Table~\ref{table_multisoliton} except those denoted as $\diamondsuit$ and $\clubsuit$ are known, and they are indicated directly below the solution kernel symbols. Solutions with $(s,s-s_z)=(i,i)$ have Petrov type X$^+$X$^-$ for some X: they have the same Petrov type on both sides and are known to admit Lorentzian sections.

\subsection{Solitonic classification of the multi-centre solutions}

The $n$-centre solutions  for $n\le 3$ studied explicitly in the previous section fit nicely into Table~\ref{table_multisoliton}. They all have Petrov type D$^+$X$^-$ for some X. The 1-centre solution is 1TN with $(s,s-s_z)=(1,2)$ in Table~\ref{table_multisoliton} and Petrov type D$^+$O$^-$, and it admits a D-to-O ALE limit to 1GH (1TN $\rightharpoondown$ 1GH) and no D-to-D ALE limit. The 2-centre solution is KN with $(s,s-s_z)=(2,2)$ in Table~\ref{table_multisoliton} and  Petrov type D$^+$D$^-$. It admits both D-to-O and D-to-D ALE limits: the D-to-O ALE limit (KN $\rightharpoondown$ 2GH) results in 2GH with Petrov type O$^+$D$^-$;  the D-to-D ALE limit (2GH $\leftharpoonup$ KN) results in 2GH with Petrov type D$^+$O$^-$. The 3-centre solution is CT with $(s,s-s_z)=(3,2)$ in Table~\ref{table_multisoliton} and Petrov type D$^+$I$^-$. It has a D-to-O ALE limit resulting in 3GH solution with Petrov type O$^+$I$^-$ (CT $\rightharpoondown$ 3GH), and a D-to-D ALE limit resulting in PD solution with Petrov type D$^+$D$^-$ (PD $\leftharpoonup$  CT).

It is of interest to understand how the $n$-centre solutions for $n>3$ fit into the above classification scheme of multi-soliton solutions on flat space presented in Table~\ref{table_multisoliton}. We will count columns from the right and use $s-s_z$ as the column index in this table. We first note that solutions in column 0 and column 1 are the multi-centre Taub--NUT and Gibbons--Hawking solutions respectively,\footnote{This statement has been verified explicitly in \cite{Chen:2015iex} for $0\le s\le 3$. A proof in the general case may be possible.} both of which have Petrov type O$^+$X$^-$ for some X. A solution $(s,s-s_z)=(i,j)$ in Table~\ref{table_multisoliton} is a special case of a solution that is directly below it with $(s,s-s_z)=(i+1,j)$: Start with an ALE solution in the table, we can follow an arrow $\nearrow$ and then an arrow $\leftharpoonup$ to reach one row up; start with an ALF solution in the table, we can follow an arrow $\leftharpoonup$ and then an arrow $\nearrow$ to reach one row up. So the Petrov type either becomes more general or stays the same when we go down each column.  A solution in Table~\ref{table_multisoliton} is a special case of the solution that sits at its bottom-left adjacent corner: we reach $(s,s-s_z)=(i,j)$ from $(s,s-s_z)=(i+1,j+1)$ by following an arrow $\nearrow$. So the Petrov type either becomes more general or stays the same when we go in the bottom-left direction in the table against the arrows $\nearrow$. These observations imply that solutions in column 4 and beyond with $4\le s-s_z\le s$ are all of Petrov type I$^+$I$^-$. 

Therefore, if the one-sided type-D Ricci-flat multi-centre solutions studied in this paper with Petrov type D$^+$X$^-$ for some X are encompassed in Table~\ref{table_multisoliton}, they can only occupy column 2 and column 3. We have seen that this is indeed true when the number of centres $n\le 3$. We expect that this is true for general $n$. With the observation that solutions in even columns  in Table~\ref{table_multisoliton} are ALF and those in odd columns are ALE, we are led to make the following conjecture:
\begin{align}
	\text{The $n$-centre solution}&\sim  (s,s-s_z)=(n,2) \text{ soliton solution ($n\ge1$)}\\
	\text{ALE $n$-centre solution}&\sim (s,s-s_z)=(n,3)	 \text{ soliton solution ($n\ge 2$)}
\end{align}
The ALE $n$-centre solution in the above conjecture refers to the D-to-D ALE limit of the $n$-centre solution. A symbol $\sim$ connecting two solutions means that they are related by some coordinate transformations and parameter redefinitions. Contained in this conjecture are the following identifications
\begin{align}
\text{The 4-centre solution}&\sim \clubsuit \label{club_conjecture}\\
 \text{ALE 4-centre solution}	&\sim \diamondsuit\label{diamond_conjecture}
\end{align}
We can prove that, in each of (\ref{club_conjecture}) and (\ref{diamond_conjecture}), the rod structures of the solutions on both sides of the symbol $\sim$ match exactly via parameter redefinitions, but we fail to map the ranges of parameters\footnote{The ranges of parameters of multi-soliton solutions on flat space were not given in \cite{Chen:2015iex}. It is a problem worth further investigation.} and prove the full equivalence of corresponding metrics. 

The above conjecture, if true, gives solitonic interpretation of the multi-centre solutions studied in this paper. In particular, it identifies the number of centres with the free-soliton number. According to the conjecture, column 2 in Table~\ref{table_multisoliton} is occupied by the $n$-centre solutions. Their D-to-D ALE limit is the decompactification limit $\leftharpoonup$ connecting column 2 and column 3. Their D-to-O ALE limit is the decompactification limit $\rightharpoondown$ connecting column 2 and column 1, producing the multi-centre Gibbons--Hawking solutions. The conjecture above, if true, completes the Petrov classification of multi-soliton solutions on flat space in table~\ref{table_multisoliton}: For $(s\ge 4, 4\le s-s_z\le s)$, all solutions have Petrov type I$^+$I$^-$; for $(s\ge4,s-s_z=2\text{ or }3)$, all solutions have Petrov type D$^+$I$^-$; for $(s\ge4,s-s_z=0\text{ or }1)$, all solutions have Petrov type O$^+$I$^-$.

\section{Discussion}

\label{sec_discussion}

We have constructed and studied one-sided type-D Ricci-flat multi-centre solutions. Many closed-form one-sided type-D Ricci-flat metrics were recovered from them explicitly. We showed that the $n$-centre solutions for $n\le 3$ occupy column 2 and column 3 in the solitonic classification table (Table~\ref{table_multisoliton}) of multi-soliton solutions on flat space, and conjectured that it is true for general $n$.

For multi-soliton solutions on flat space classified in Table~\ref{table_multisoliton}, it was observed in \cite{Chen:2015iex} that a solution with $(s,s-s_z)=(i+2,j+2)$ is a two-soliton transformation on a solution with $(s,s-s_z)=(i,j)$: PD is a two-soliton transformation on 1GH, CT is a two-soliton transformation on 1TN, etc.. This observation can be proved using the equivalence between adding solitons at one time and adding the same solitons successively in ISM constructions. Therefore, each solution in column 2 in Table~\ref{table_multisoliton} can be obtained from a solution in column 0 by a two-soliton transformation, and each solution in column 3 can be obtained from a solution in column 1 by a two-soliton transformation. The conjecture we made in the previous section can be proved if we can prove the following: (a) multi-soliton solutions on flat space with $s-s_z=0,1$ are one-sided type O and (b) a two-soliton transformation changes a one-sided type-O solution to a one-sided type-D solution.

The conjecture we made, if true, has implications in the search of new gravitational instantons with four or more turning points from Table~\ref{table_multisoliton}. In the case of four turning points, we should search for them in the fourth row, with $s=4$, in this table. 4TN and 4GH are well understood to describe gravitational instantons if all the centre charges are equal. The solutions \textbf{dKN}, $\diamondsuit$ and $\clubsuit$ with solitonic representations described as
\begin{align}
	\text{\bf dKN}&\sim[8-2-2]\sim [4-0-0],\\
	\diamondsuit&\sim[7-2-1]\sim [5-1-0],\\
	\clubsuit&\sim [8-3-1]\sim[6-2-0],	
\end{align}
respectively, are less understood. If our conjecture holds, both the solution $\diamondsuit$ and the solution $\clubsuit$ are one-sided type D. The classification of Biquard and Gauduchon \cite{Biquard:2021gwj}, or the result proved in Appendix \ref{sec_proof_nge3}, then implies that neither the solution $\diamondsuit$ nor the solution $\clubsuit$ contains gravitational instantons.\footnote{With ranges of parameters of the two sides of each of (\ref{club_conjecture}) and (\ref{diamond_conjecture}) not mapped, the equivalence of the rod structures of the two sides that we mentioned earlier (see comments below (\ref{diamond_conjecture})), together with the result proved in Appendix~\ref{sec_proof_nge3}, is not sufficient to exclude the possibility that the solutions $\diamondsuit$ and $\clubsuit$ might contain gravitational instantons. The stronger result of equivalence of metrics on both sides implied by the conjecture is needed.}\textsuperscript{,}\footnote{This does not exclude the possibility that the solutions $\diamondsuit$ and $\clubsuit$ can give rise to gravitational instantons once charged say with Maxwell fields.} So the most promising candidate in Table~\ref{table_multisoliton} to extract 4-turning-point gravitational instantons from is the double Kerr--NUT solution with $(s,s-s_z)=(4,4)$. {A new AF gravitational instanton with four turning points whose existence was proved by Li and Sun \cite{Li:2025} has recently been constructed explicitly by Teo \cite{Teo:2026kej}.}

We note that in Table~\ref{table_multisoliton}, for $s\le 3$ (with 3TN excluded), the Petrov type and the asymptotic structure (being ALF/ALE) together determine a unique solution. This is no longer true for $s\ge 4$. In particular, all solutions with $s\ge s-s_z\ge 4$ are algebraically general, with Petrov type I$^+$I$^-$. It is natural to wonder if a finer algebraic classification of the Weyl tensor, perhaps one that can be applied to submanifolds like $S^2$-cycles corresponding to rods in the rod structure, is possible. Such a classification would be helpful to characterise multi-soliton solutions with large soliton numbers. Uniqueness and stability of multi-soliton solutions in Table~\ref{table_multisoliton} may be considered. The harmonic map formulation of toric Ricci-flat metrics \cite{Kunduri:2021xiv,Li:2025} may be useful to prove (non-)uniqueness of these solutions.

A promising extension of our work is to construct one-sided type-D multi-centre solutions with Maxwell fields. As a by-product in our derivation, we have obtained scalar-flat conformally K{\"a}hler metrics (see Eqs.~(\ref{metric_Lebrun_Tod}), (\ref{Kahler_form1}), (\ref{Toda_3D}), (\ref{monopolo_equation}) and (\ref{Equation_for_A2})) and the trace-free Ricci curvature operator (\ref{trace_free_Ricci}) of these metrics. Our calculations agree with earlier results performed by Araneda \cite{Araneda:2023xnv} and Tod \cite{Tod:2024msz}. It is known that these metrics can be turned into solutions of Einstein--Maxwell equations \cite{Flaherty:1978wh,Araneda:2023xnv}: the self-dual Maxwell field is essentially given by the fundamental two-form (rescaled to be closed by an appropriate conformal factor) and the anti-self-dual Maxwell field is essentially given by the Ricci two-form (rescaled to be closed). The function $W$ is no longer constrained by (\ref{W_equation1}).  Instead, it is now only subject to the monopole equation (\ref{monopolo_equation}), which can be reduced to another harmonic potential (distinct from the generating potential $x$) in the auxiliary 3D flat space (\ref{3D_auxiliary_space1}), see \cite{Tod:2024msz}.  One-sided type-D multi-centre solutions charged with Maxwell fields may be obtained from a certain choice of this second harmonic potential. The charged 2-centre solution should reproduce the Kerr--Newman solution (with NUT charge). The charged 3-centre solution should reproduce or generalise the solution found by Araneda and Dunajski \cite{Araneda:2025lak} with a charged new AF instanton emerging as a special case.

Calderbank and Pedersen \cite{Calderbank:2001uz} obtained the generalisation of the multi-centre Gibbons--Hawking metrics with a cosmological constant, which can be written in the following form as given in \cite{Casteill:2001zk}:
\begin{align}
	\label{metric_CP}
	\dif s^2&=w^{-2}[W^{-1}(\dif\tau+A)^2+W\rho^2\dif\phi^2+WL(\dif\rho^2+\dif z^2)],\\
	L&=w_\rho^2+w_z^2,\quad W=1-\frac{ ww_\rho}{\rho L},\quad A=\left(z-\frac{ww_z}{L}\right)\dif\phi.\nonumber
\end{align} 
The function $w$ satisfies (\ref{w_harmonicity}) and, for multi-centre solutions,  it is given by  (\ref{w_of_multisoliton}). These solutions are anti-self-dual (or self-dual)  and satisfy $R_{\mu\nu}+3g_{\mu\nu}=0$. The metrics (\ref{metric_onesided_typeD}) and (\ref{metric_CP}) are similar in structure. The generalisation of these metrics to one-sided type-D multi-centre metrics with a cosmological constant may be worth pursuing. Tod \cite{Tod:2020ual} pointed out that when both the one-sided type-D condition and a cosmological constant are imposed, the field equations are not integrable.

It is well known that the multi-centre Gibbons--Hawking metrics admit generalisations in which the centres are not all aligned along any axis, and the metrics in this case have only one Killing vector. Do one-sided type-D Ricci-flat multi-centre metrics have generalisations with centres not all aligned on any axis? There is also only one Killing vector, guaranteed by the one-sided type-D condition. The Weyl--Papapetrou coordinates are of course inadequate to describe such metrics (if they exist at all).

\section*{Acknowledgements}

The author would like to thank the Simons Center for Geometry and Physics, Stony Brook University, for support during the program Einstein 4-manifolds and gravitational instantons, where part of the research for this paper was performed. He would like to thank the organizers of this program. His participation in this program was also supported by PAP General Purpose Fund C230069300 from Nanyang Technological University. He wishes to thank Paul Tod for helpful discussions on one-sided type-D Ricci-flat metrics. {He would like to thank Edward Teo for comments on the manuscript and wishes to thank the anonymous referee for various suggestions that helped improve the presentation of the paper.}

\appendix

\section{Proof of  $n\le 3$ for gravitational instantons}

\label{sec_proof_nge3}

For the $n$-centre solution to describe gravitational instantons, it is necessary that $n\le 3$. This result was first proved by Biquard and Gauduchon \cite{Biquard:2021gwj}. In this appendix, we give a simplified proof.

For a toric Ricci-flat solution to describe gravitational instantons, any ordered pair of adjacent rod directions in its rod structure must be related by $GL(2,\mathbb{Z})$ transformations \cite{Chen:2010zu} to the pair of independent generators of the toric isometry. Then Eq.~(\ref{M_K}) implies that, in the case of the $n$-centre solution, these $GL(2,\mathbb{Z})$ transformations are actually $SL(2,\mathbb{Z})$ transformations. Thus, the matrices $M_k$, given by (\ref{M_K}), must be congruent to $SL(2,\mathbb{Z})$ matrices:
\begin{equation}
	NM_k N^T\in SL(2,\mathbb{Z}),
\end{equation}
for all $1\le k\le n$ with a matrix $N$ independent of $k$. It follows that $\det M_k$ is $k$-independent:
\begin{equation}
	\label{determinant_of_Mk}
	\det M_k=(\det N)^{-2},
\end{equation}
for all $1\le k\le n$. Alternatively, we can consider neighbouring matrix pairs $(M_{k-1},M_k)$ and impose the condition $M_{k}=M_{k-1}S_k$ for some $k$-dependent matrix $S_k\in SL(2,\mathbb{Z})$. Both $M_{k-1}$ and $M_k$ contain $v_k$ as a column. The form of $S_k$ can be restricted:
\begin{align}
	\label{SL2Z}
	M_{k}=M_{k-1}S_k: 	(v_k,v_{k+1})=(v_{k-1},v_k)S_k=(v_{k-1},v_k)\begin{pmatrix}
		0	& q_k\\1&p_k
	\end{pmatrix},
\end{align}
where $p_k$ and $q_k$ are two integers. The condition $S_k\in SL(2,\mathbb{Z})$ implies that $q_k=-1$. Then (\ref{SL2Z}) is equivalent to
\begin{empheq}[left={v_{k+1}+v_{k-1}=p_kv_k \quad\Leftrightarrow\quad\empheqlbrace}]{align}
	& b_{k-1}+b_{k+1}=p_kb_k  \label{instanton_condition1} \\
	& a_{k-1}+a_{k+1}=p_ka_k  \label{instanton_condition2}
\end{empheq}
where $a_k$ and $b_k$ are defined in (\ref{rod_directions}) and  $p_k\in \mathbb{Z}$. If we eliminate $p_k$ by combining (\ref{instanton_condition1}) and (\ref{instanton_condition2}), we obtain the following equation 
\begin{equation}
		\label{M_equation}
	\det M_{k-1}=\det M_k,
\end{equation}
which also follows from (\ref{determinant_of_Mk}) directly.

Imposing the two conditions (\ref{instanton_condition1}) and (\ref{instanton_condition2}) on the multi-centre solutions leads to Biquard--Gauduchon classification of Hermitian and non-K{\"a}hler gravitational instantons. A key result of their classification is that $n\le3$ for the multi-centre solutions to describe gravitational instantons. We give a proof of this result here. We use the necessary conditions (\ref{instanton_condition1}) and (\ref{M_equation}). In view of (\ref{M_K}), Eq.~(\ref{M_equation}) implies that
\begin{equation}
	\sqrt{\alpha_{k-1}}\left[\beta+\kminusonesum\alpha_i (z_{k-1}-z_i)\right]=\sqrt{\alpha_k}\left[\beta+\ksum\alpha_i(z_k-z_i)\right],
\end{equation}
which can be rearranged to give
\begin{equation}
	\frac{\sqrt{\alpha_k}}{\sqrt{\alpha_{k-1}}}=1-\frac{(z_k-z_{k-1})\ksum\alpha_i}{\beta+\ksum\alpha_i(z_k-z_i)}.
\end{equation}
Now use (\ref{jump_in_b}) to express the above equation as
\begin{equation}
	\sqrt{\frac{b_{k+1}-b_k}{b_k-b_{k-1}}}=1-\frac{(z_{k}-z_{k-1})k_2^{-1}b_k}{\beta+\ksum\alpha_i(z_k-z_i)}= 1-\gamma b_k, \quad \gamma\equiv\frac{(z_{k}-z_{k-1})k_2^{-1}}{\beta+\ksum\alpha_i(z_k-z_i)}.
\end{equation}
We notice that $\gamma>0$ as both its numerator and denominator are positive. By substituting (\ref{instanton_condition1}) in the above equation, we obtain
\begin{equation}
	\label{instanton_condition3}
	\sqrt{1+\frac{(p_k-2)b_k}{b_k-b_{k-1}}}=1-\gamma b_k.
\end{equation}
This equation contradicts (\ref{instanton_condition1}) if $0\le b_{k-1}<b_k<b_{k+1}$ or $b_{k-1}<b_k<b_{k+1}\le 0$. In both cases, (\ref{instanton_condition1}) implies that $p_k\ge 2$. If $0\le b_{k-1}<b_k<b_{k+1}$, then the left-hand side of (\ref{instanton_condition3}) is greater than 1, contradicting the right-hand side which is less than 1. If $b_{k-1}<b_k<b_{k+1}\le 0$, the left-hand side is less than 1, contradicting again the right-hand side which is now greater than 1. Combining these two cases, we thus deduce that $b_{k-1}$ and $b_{k+1}$ must have opposite signs. Therefore, $b_3>0$ in view of $b_1<0$, and $b_{n-1}<0$ in view of $b_{n+1}>0$. So we have obtained the inequality
\begin{equation}
	b_{n-1}<0<b_3.
\end{equation}
It then follows that $n-1\le 2$, and therefore $n\le 3$.

\bigskip\bigskip

{\renewcommand{\Large}{\large}
}


\begin{thebibliography}{99}

\bibitem{Gibbons:1976ue}
G.~W.~Gibbons and S.~W.~Hawking,
``Action integrals and partition functions in quantum gravity,''
\href{https://doi.org/10.1103/PhysRevD.15.2752}{Phys. Rev. D \textbf{15} (1977), 2752-2756.}


\bibitem{Hawking:1976jb}
S.~W.~Hawking,
``Gravitational instantons,''
\href{https://doi.org/10.1016/0375-9601(77)90386-3}{Phys. Lett. A \textbf{60} (1977), 81.}

\bibitem{Przanowski:1983ybm}
M.~Przanowski and B.~Broda,
``Locally K{\"a}hler gravitational instantons,''
Acta Phys. Polon. B \textbf{14} (1983), 637-661.

\bibitem{Derdzinski}
A. Derdzi\'{n}ski, ``Self-dual K{\"a}hler manifolds and Einstein manifolds of dimension four,'' Compos. Math. \textbf{49} (1983), 405–433.

\bibitem{Goldblatt:1994rx}
E.~Goldblatt,
``A Newman-Penrose formalism for gravitational instantons,'' \href{https://doi.org/10.1007/BF02106666}{
Gen. Rel. Grav. \textbf{26} (1994), 979-997.}


\bibitem{Gibbons:1978tef}
G.~W.~Gibbons and S.~W.~Hawking,
``Gravitational multi-instantons,''
 \href{https://doi.org/10.1016/0370-2693(78)90478-1}{Phys. Lett. B \textbf{78} (1978), 430.}



\bibitem{Newman:1963yy}
E.~Newman, L.~Tamburino and T.~Unti,
``Empty space generalization of the Schwarzschild metric,''
\href{https://doi.org/10.1063/1.1704018}{J. Math. Phys. \textbf{4} (1963), 915.}



\bibitem{Eguchi:1978xp}
T.~Eguchi and A.~J.~Hanson,
``Asymptotically flat selfdual solutions to Euclidean gravity,''
\href{https://doi.org/10.1016/0370-2693(78)90566-X}{Phys. Lett. B \textbf{74} (1978), 249-251.}



\bibitem{Page:1978hdy}
D.~N.~Page,
``Taub-NUT instanton with an horizon,''
\href{https://doi.org/10.1016/0370-2693(78)90016-3}{Phys. Lett. B \textbf{78} (1978), 249-251.}


\bibitem{Chen:2011tc}
Y.~Chen and E.~Teo,
``A new AF gravitational instanton,''
\href{https://doi.org/10.1016/j.physletb.2011.07.076}{Phys. Lett. B \textbf{703} (2011), 359-362}
[arXiv:1107.0763 [gr-qc]].


\bibitem{Aksteiner:2021fae}
S.~Aksteiner and L.~Andersson,
``Gravitational instantons and special geometry,''
\href{https://doi.org/10.4310/jdg/1729092451}{J. Diff. Geom. \textbf{128} (2024) no.3, 928-958}
[arXiv:2112.11863 [gr-qc]].


\bibitem{Biquard:2021gwj}
O.~Biquard and P.~Gauduchon,
``On toric Hermitian ALF gravitational instantons,''
\href{https://doi.org/10.1007/s00220-022-04562-z}{Commun. Math. Phys. \textbf{399} (2023) no.1, 389-422}
[arXiv:2112.12711 [math.DG]].

\bibitem{Carter:1968ks}
B.~Carter,
``Hamilton-Jacobi and Schrodinger separable solutions of Einstein's equations,''
\href{https://doi.org/10.1007/BF03399503}{Commun. Math. Phys. \textbf{10} (1968) no.4, 280-310.}


\bibitem{Gibbons:1979nf}
G.~W.~Gibbons and M.~J.~Perry,
``New gravitational instantons and their interactions,''
\href{https://doi.org/10.1103/PhysRevD.22.313}{Phys. Rev. D \textbf{22} (1980), 313.}


\bibitem{Chen:2015vva}
Y.~Chen and E.~Teo,
``Five-parameter class of solutions to the vacuum Einstein equations,''
\href{https://doi.org/10.1103/PhysRevD.91.124005}{Phys. Rev. D \textbf{91} (2015) no.12, 124005}
[arXiv:1504.01235 [gr-qc]].


\bibitem{Przanowski:1984}
M.~Przanowski,
``One-sided type-D gravitational instantons,''
\href{https://doi.org/10.1007/BF00762933}{Gen. Rel. Grav. \textbf{16} (1984), 797-803.}

\bibitem{Tod:2020ual}
P.~Tod,
``One-sided type-D Ricci-flat metrics,'' N.~Y.~J.~Math., \textbf{32} (2026), 282-302 [arXiv:2003.03234 [gr-qc]].


\bibitem{Lebrun:1991}
C.~LeBrun,
``Explicit self-dual metrics on $\mathbb{C}P_2\#...\#\mathbb{C}P_2$,''
\href{https://doi.org/10.4310/jdg/1214446999}{J. Diff. Geom. \textbf{34} (1991) no.1, 223-253.}

\bibitem{Biquard:2023rbr}
O.~Biquard and P.~Gauduchon,
``About a family of ALF instantons with conical singularities,''
\href{https://doi.org/10.3842/SIGMA.2023.079}{SIGMA \textbf{19} (2023), 079}
[arXiv:2306.11110 [math.DG]].


\bibitem{Tod:2024msz}
P.~Tod,
``One-sided type-D metrics with aligned Einstein-Maxwell,'' N.~Y.~J.~Math., \textbf{32} (2026), 303-320
[arXiv:2410.13410 [gr-qc]].


\bibitem{Penrose:1985bww}
R.~Penrose and W.~Rindler,
``Spinors and space-time. Vol. 1: Two-spinor calculus and relativistic fields,''
Cambridge University Press, 1984.

\bibitem{Penrose:1986ca}
R.~Penrose and W.~Rindler,
``Spinors and spacetime. Vol. 2: Spinor and twistor methods in space-time geometry,''
Cambridge University Press, 1988.

\bibitem{Dunajski:2010zz}
M.~Dunajski,
``Solitons, instantons, and twistors,'' Oxford University Press, 2010.



\bibitem{Araneda:2025uqo}
B.~Araneda and J.~Lucietti,
``All toric Hermitian ALE gravitational instantons,''
[arXiv:2510.09291 [math.DG]].


\bibitem{Belinski:2001ph}
V.~Belinski and E.~Verdaguer,
``Gravitational solitons,''
Cambridge University Press, 2005.

\bibitem{Belinsky:1971nt}
V.~A. Belinsky and V.~E. Zakharov, { {``Integration of the Einstein equations
by the inverse scattering problem technique and the calculation of the exact
soliton solutions,''}}  {Sov. Phys. JETP} {\bf 48} (1978) 985-994.


\bibitem{Chen:2015iex}
Y.~Chen,
``Gravitational multisoliton solutions on flat space,''
\href{https://doi.org/10.1103/PhysRevD.93.044021}{Phys. Rev. D \textbf{93} (2016) no.4, 044021}
[arXiv:1512.00032 [gr-qc]].



\bibitem{Emparan:2001wk}
R.~Emparan and H.~S.~Reall,
``Generalized Weyl solutions,''
\href{https://doi.org/10.1103/PhysRevD.65.084025}{Phys. Rev. D \textbf{65} (2002), 084025}
[arXiv:hep-th/0110258 [hep-th]].


\bibitem{Harmark:2004rm}
T.~Harmark,
``Stationary and axisymmetric solutions of higher-dimensional general relativity,''
\href{https://doi.org/10.1103/PhysRevD.70.124002}{Phys. Rev. D \textbf{70} (2004), 124002}
[arXiv:hep-th/0408141 [hep-th]].

\bibitem{Chen:2010zu}
Y.~Chen and E.~Teo,
``Rod-structure classification of gravitational instantons with $U(1)\times U(1)$ isometry,''
\href{https://doi.org/10.1016/j.nuclphysb.2010.05.017}{Nucl. Phys. B \textbf{838} (2010), 207-237}
[arXiv:1004.2750 [gr-qc]].


\bibitem{Plebanski:1976gy}
J.~F.~Pleba\'nski and M.~Demia\'nski,
``Rotating, charged, and uniformly accelerating mass in general relativity,''
\href{https://doi.org/10.1016/0003-4916(76)90240-2}{Annals Phys. \textbf{98} (1976), 98-127.}




\bibitem{Li:2025}
M.~Li and S.~Sun,
``Gravitational instantons and harmonic maps,''
[arXiv:2507.15284 [math.DG]].


\bibitem{Teo:2026kej}
{E.~Teo,
``New asymptotically flat gravitational instanton,''
[arXiv:2605.28542 [gr-qc]].}


\bibitem{Kunduri:2021xiv}
H.~K.~Kunduri and J.~Lucietti,
``Existence and uniqueness of asymptotically flat toric gravitational instantons,''
\href{https://doi.org/10.1007/s11005-021-01475-1}{Lett. Math. Phys. \textbf{111} (2021) no.5, 133}
[arXiv:2107.02540 [math.DG]].


\bibitem{Araneda:2023xnv}
B.~Araneda,
``Hidden symmetries of generalised gravitational instantons,''
\href{https://doi.org/10.1007/s00023-024-01515-1}{Ann.\ Henri Poincar\'e \textbf{26} (2025) no.11, 4021-4049}
[arXiv:2309.05617 [gr-qc]].


\bibitem{Flaherty:1978wh}
E.~J.~Flaherty,
``The nonlinear graviton in interaction with a photon,''
\href{https://doi.org/10.1007/BF00784657}{Gen. Rel. Grav. \textbf{9} (1978), 961-978.}


\bibitem{Araneda:2025lak}
B.~Araneda and M.~Dunajski,
``New asymptotically flat Einstein-Maxwell instantons,''
\href{https://doi.org/10.1103/f3ls-znl6}{Phys. Rev. Lett. \textbf{135} (2025) no.24, 241501}
[arXiv:2510.06458 [gr-qc]].

\bibitem{Calderbank:2001uz}
D.~M.~J.~Calderbank and H.~Pedersen,
``Selfdual Einstein metrics with torus symmetry,''
\href{https://doi.org/10.4310/jdg/1090351125}{J. Diff. Geom. \textbf{60} (2002) no.3, 485-521}
[arXiv:math/0105263 [math.DG]].

\bibitem{Casteill:2001zk}
P.~Y.~Casteill, E.~Ivanov and G.~Valent,
``$U(1)\times U(1)$ quaternionic metrics from harmonic superspace,''
\href{https://doi.org/10.1016/S0550-3213(02)00013-5}{Nucl. Phys. B \textbf{627} (2002), 403-444}
[arXiv:hep-th/0110280 [hep-th]].


\end{thebibliography}
\end{document}